\documentclass[12pt]{article}
\usepackage{graphicx}
\usepackage{amssymb}

\newcommand{\sect}[1]{\section{#1}\setcounter{equation}{0}}

\font\mbn=msbm10 scaled \magstep1
\font\mbs=msbm7 scaled \magstep1
\font\mbss=msbm5 scaled \magstep1
\newfam\mbff
\textfont\mbff=\mbn
\scriptfont\mbff=\mbs
\scriptscriptfont\mbff=\mbss\def\mbf{\fam\mbff}
\def\Re{{\mbf R}}

\def\Z{{\mbf Z}}
\def\Co{{\mbf C}}

\def\Di{{\mbf D}}

\def\N{{\mbf N}}
\def\H{{\mbf H}}
\newtheorem{Th}{Theorem}[section]
\newtheorem{Lm}[Th]{Lemma}
\newtheorem{C}[Th]{Corollary}

\newtheorem{D}[Th]{Definition}

\newtheorem{R}[Th]{Remark}
\newtheorem{E}[Th]{Example}

\author{Alexander Brudnyi\thanks{Research supported in part by NSERC.
}\\
Department of Mathematics and Statistics\\
University of Calgary, Calgary\\
Canada\\
\\
Damir Kinzebulatov\\
Department of Mathematics and Statistics\\
University of Calgary, Calgary\\
Canada}
\title{On Uniform Subalgebras of $L^{\infty}$ on the Unit Circle Generated by Almost Periodic Functions}
\date{}
\begin{document}
\maketitle
\begin{abstract}
{In the present paper we introduce analogs of almost periodic functions for the unit circle. We study certain uniform algebras generated by such functions, prove corona theorems for them and describe their maximal ideal spaces.}
\end{abstract}
\sect{\hspace*{-1em}. Formulation of Main Results}
{\bf 1.1.} The classical almost periodic functions on the real line as first introduced by \penalty-10000 H. Bohr in the 1920s play an important role in various areas of Analysis. In the present paper we define analogs of almost periodic functions on the unit circle. We study certain uniform algebras generated by such functions. In particular, we describe in these terms some uniform subalgebras of the algebra $H^{\infty}$ of bounded holomorphic functions on the open unit disk $\Di\subset\Co$, having, in a sense, the weakest possible discontinuities on the boundary $\partial\Di$.

To formulate the main results of the paper we first recall the definition of almost periodic functions, see [B].
\begin{D}\label{d1}
A continuous function $f:\Re\to\Co$ is called almost periodic if, for any $\epsilon>0$, there exists $l(\epsilon)>0$  such that for every $t_0 \in \Re$ the interval $[t_{0},t_{0}+l(\epsilon)]$ contains at least one number $\tau$  for which
$$
|f(t)-f(t+\tau)|<\epsilon\ \ \ {\rm for\ all}\ \ \ t\in\Re.
$$
\end{D}
It is well known that every almost periodic function $f$ is uniformly continuous and is the uniform limit of a sequence of exponential polynomials
$\{q_{n}\}_{n\in\N}$ where $q_{n}(t):=\sum_{k=1}^{n}c_{kn}e^{i\lambda_{kn}t}$, $c_{kn}\in\Co$, $\lambda_{kn}\in\Re$, $1\leq k\leq n$, and  $i:=\sqrt{-1}$.

In what follows we consider $\partial\Di$ with the {\em counterclockwise orientation}.
For $t_0 \in \Re$ let $\gamma_{t_{0}^{k}}(s):=\{e^{i(t_{0}+kt)}\ :\ 0<t<s\leq 2\pi\}\subset \partial\Di$,
$k\in\{-1,1\}$, be two open arcs having $e^{it_{0}}$ as the right or the left
endpoints with respect to the chosen orientation, respectively.

Let us define almost periodic functions on open arcs of $\partial\Di$.
\begin{D}\label{d2}
A continuous function $f_{k}: \gamma_{t_{0}^{k}}(s)\to\Co$, $k\in\{-1,1\}$, is said to be
almost periodic if the function $\widehat f_{k}:(-\infty, 0)\to\Co$,\
$\widehat f_k(t):=f_k(e^{i(t_{0}+kse^{t})})$, is the restriction of an almost
periodic function on $\Re$.
\end{D}
\begin{E}\label{ex1}
{\rm The function $e^{i\lambda\log_{t_{0}^{k}}}$, $\lambda\in\Re$,
where
$$
\log_{t_{0}^{k}}(e^{i(t_{0}+kt)}):=\ln t,\ \ \ 0<t<2\pi,\ \ k\in\{-1,1\},
$$
is almost periodic in the sense of this definition on 
$\gamma_{t_{0}^{1}}(2\pi)= \gamma_{t_{0}^{-1}}(2\pi)$.}
\end{E}

By $AP(\partial\Di)\subset L^{\infty}(\partial\Di)$ we denote the uniform subalgebra of functions $f$ such that for each $t_{0}$ and any $\epsilon>0$ there are a number $s:=s(t_{0},\epsilon)\in (0,\pi)$ and almost periodic functions $f_{k}:\gamma_{t_{0}^{k}}(s)\to\Co$, $k\in\{-1,1\}$, such that
\begin{equation}\label{e1}
{\rm ess}\!\!\!\!\!\sup_{z\in \gamma_{t_{0}^{1}}(s)}|f(z)-f_{1}(z)|<\epsilon\ \ \ {\rm and}
\ \ \ {\rm ess}\!\!\!\!\!\sup_{z \in \gamma_{t_{0}^{-1}}(s)}|f(z)-f_{-1}(z)|<\epsilon.
\end{equation}

Let $S\subset \partial\Di$ be a nonempty closed subset. By $AP(S)\subset AP(\partial\Di)$ we denote the uniform algebra of functions from $AP(\partial\Di)$ continuous on $\partial\Di\setminus S$. 

Fix a real continuous
function $g$ on $\gamma_{t_{0}^{1}}(2\pi)$ such that 
$$
\lim_{t\to 0+}g(e^{it})=1,\ \ \ \ \ \lim_{t\to 2\pi -}g(e^{it})=0 
$$
and $g(e^{it})$ is decreasing for $0<t<2\pi$. We set $g_{t_{0}}(e^{it}):=g(e^{i(t_{0}+t)})$.
\begin{Th}\label{te1}
The algebra $AP(S)$ is the uniform closure in $L^{\infty}(\partial\Di)$ of the
algebra of complex polynomials in variables $g_{t_{0}}$ and
$e^{i\lambda\log_{t_{0}^{k}}}$, $\lambda\in\Re$, $e^{it_{0}}\in S$, $k\in\{-1,1\}$.
\end{Th}

Let $\phi:\partial\Di\to\partial\Di$ be a $C^{1}$ diffeomorphism. By $\phi^{*}:C(\partial\Di)\to C(\partial\Di)$, 
$\phi^{*}(f):=f\circ\phi$, we denote 
the pullback by $\phi$. Set $\widetilde S:=\phi(S)$. As a consequence of Theorem \ref{te1} we obtain
\begin{C}\label{cor1}
$\phi^{*}$ maps $AP(\widetilde S)$ isomorphically onto $AP(S)$.
\end{C}  
{\bf 1.2.} We say that a complex-valued function $g\in L^{\infty}(\partial\Di)$ has a discontinuity of the {\em first kind} at $x_{0}$, if the one-sided limits of $g$ at $x_{0}$ exist but have distinct values. 
For a closed subset $S\subset \partial\Di$ by $R_{S}\subset L^{\infty}(\partial\Di)$ we denote the uniform algebra of complex functions allowing discontinuities of the first kind at points of $S$ and continuous on $\partial\Di\setminus S$. Elements from $R_{S}$ are often referred to as {\em regulated functions} [D].  Clearly, $R_{S}\hookrightarrow AP(S)$. Also, we will show (see Lemma \ref{unif1} below) that $R_{S}$ is the uniform closure of the algebra generated by all possible subalgebras $R_{F}$ with finite $F\subset S$.

Let ${\cal M}(AP(S))$ be the {\em  maximal ideal space} of $AP(S)$, that is, the space of all {\em characters} (=\ nonzero homomorphisms $AP(S)\to \Co$) on $AP(S)$ equipped with the weak$^{*}$-topology (also known as the {\em Gelfand topology}) inherited from $(AP(S))^{*}$. By definition, ${\cal M}(AP(S))$ is a compact Hausdorff space. The main result formulated in this section describes the topological structure of ${\cal M}(AP(S))$.

Let us consider continuous embeddings of uniform algebras
$$
C(\partial\Di)\hookrightarrow R_{S}\hookrightarrow AP(S).
$$
The dual maps to these embeddings determine continuous surjective maps of the
corresponding maximal ideal spaces:
$$
{\cal M}(AP(S))\stackrel{r_{S}}{\longrightarrow}{\cal M}(R_{S})
\stackrel{c_{S}}{\longrightarrow}{\cal M}(C(\partial\Di))\cong \partial\Di.
$$
\begin{Th}\label{te2}
\begin{itemize}
\item[(1)]
For each $z\in S$ preimage $c_{S}^{-1}(z)$ consists of two points $z_{+}$ and
$z_{-}$ which are naturally identified with counterclockwise and clockwise orientations of $\partial\Di$ at $z$.
\item[(2)]
$c_{S}:{\cal M}(R_{S})\setminus c_{S}^{-1}(S)\to \partial\Di\setminus S$ is a homeomorphism.

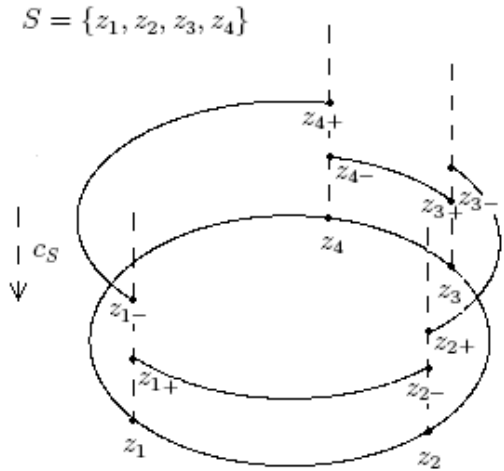
\begin{center}
\includegraphics[scale=0.7]{im1.eps}
\end{center}
{\rm Figure 1.} For a given $S=\{z_1,z_2,z_3,z_4\}$ we have the homeomorphism $c_S$ of 
\begin{displaymath}
{\cal M}(R_S) \setminus \{z_{1+},z_{1-},z_{2+},z_{2-},z_{3+},z_{3-},z_{4+},z_{4-}\}
\end{displaymath}
and $\partial \Di \setminus \{z_1,z_2,z_3,z_4\}$ where $c_S(z_{i+})=c_S(z_{i-})=z_i$.

\item[(3)]
$c_{S}^{-1}(S)\subset {\cal M}(R_{S})$ is a totally disconnected compact Hausdorff space.
\item[(4)]
For each $\xi\in c_{S}^{-1}(S)$ preimage $r_{S}^{-1}(\xi)$ is homeomorphic to the Bohr compactification $b\Re$ of $\Re$.
\item[(5)]
The map $r_{S}: {\cal M}(AP(S))\setminus (c_{S}\circ r_{S})^{-1}(S)\to {\cal M}(R_{S})\setminus c_{S}^{-1}(S)\cong \partial\Di\setminus S$ is a homeomorphism.
\end{itemize}
\end{Th}

\begin{center}
\includegraphics[scale=0.7]{im2.eps}
\end{center}
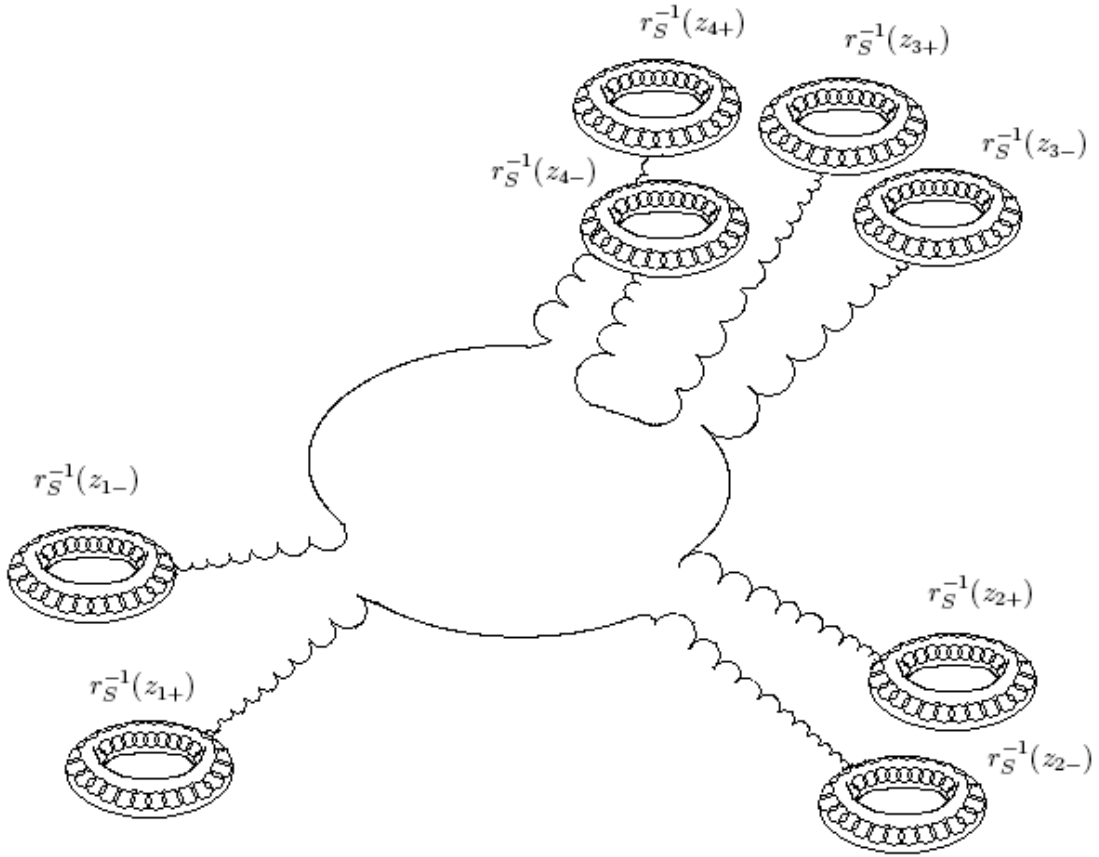
\textit{{\rm Figure 2.} Given an $n$ point set $S$, the maximal ideal space ${\cal M}(AP(S))$ is the union of $\partial \Di \setminus S$ and $2n$ Bohr compactifications of $\Re$ that can be viewed as infinite dimensional "tori", where the spirals joining the arcs and the tori are meant to indicate
(in a figurative manner) that there is an influence of the topology of the
Bohr compactifications on the topology of the arcs.} \\

Let us recall that $b\Re$ is a compact abelian topological group homeomorphic to the maximal ideal space of the algebra of continuous almost periodic functions on $\Re$. Also, it follows straightforwardly from (2)-(5) that 
\begin{itemize}
\item[(6)]
The covering dimension of ${\cal M}(AP(S))$ is $\infty$.
\item[(7)]
For a continuous map $\phi:T\to (c_{S}\circ r_{S})^{-1}(S)$ of a connected topological space $T$, there is a point $\xi\in c_{S}^{-1}(S)$ such that
$\phi(T)\subset r_{S}^{-1}(\xi)$.
\end{itemize}
{\bf 1.3.} 
By $A_{0}\subset H^{\infty}$ we denote the {\em disk-algebra}, i.e., the algebra of functions continuous on the closure $\overline{\Di}$ and holomorphic in $\Di$. Also,  by $f|_{\partial\Di}$ we denote the boundary values of $f\in C(\Di)$ (in case they exist). In the present part we describe uniform subalgebras of $H^{\infty}$ generated by almost periodic functions.
These subalgebras contain $A_{0}$ and have, in a sense, the weakest possible discontinuities on $\partial\Di$.

Suppose that $S$ contains at least $2$ points. By $A_{S}\subset H^{\infty}$ we denote the uniform closure of the algebra generated by $A_{0}$ and by holomorphic functions of the form $e^{f}$ where ${\rm Re}f|_{\partial\Di}$ is a finite linear combination with real coefficients of characteristic functions of closed arcs whose endpoints belong to $S$. If $S$ consists of a single point we determine $A_{S}\subset H^{\infty}$ to be the uniform closure of the algebra generated by $A_{0}$ and functions $ge^{\lambda f}$, $\lambda\in\Re$, where ${\rm Re}f|_{\partial\Di}$ is the characteristic function of a closed arc with an endpoint at $S$ and $g\in A_{0}$ is a function  such that $ge^{f}$ has discontinuity on $S$ only.
In the following result we naturally identify $A_{S}$ and $H^{\infty}$ with the algebras of their boundary values. 
\begin{Th}\label{te3}
$$
A_{S}=AP(S)\cap H^{\infty}.
$$
\end{Th}
\begin{R}\label{log}
{\rm Suppose that $F\subset\partial\Di$ contains at least 2 points. Let $e^{\lambda f}\in A_{S}$, $\lambda\in\Re$, where ${\rm Re}\ \! f$ is the characteristic function of an arc $[x,y]$ with $x,y\in S$. Let
$\phi_{x,y}:\Di\to\H_{+}$ be the bilinear map onto the upper half-plane that maps $x$ to $0$, the midpoint of the arc $[x,y]$ to $1$ and $y$ to $\infty$. Then
there is a constant $C$ such that
$$
e^{\lambda f(z)}=e^{-\frac{i\lambda}{\pi}{\rm Log}\ \! \phi_{x,y}(z)+\lambda C},\ \ \
z\in\Di,
$$
where ${\rm Log}$ denotes the principal branch of the logarithmic function. Thus from Theorem \ref{te3} it follows that the algebra $AP(S)\cap H^{\infty}$ is the uniform closure of the algebra generated by $A_{0}$ and the functions
$e^{i\lambda ({\rm Log}\circ\phi_{x,y})}$, $\lambda\in\Re$, $x,y\in S$.}
\end{R}

The following example shows that if $S$ is an infinite set, $A_{S}$ does not coincide with the algebra generated by functions $e^{f}$ where ${\rm Re}f\in R_{S}$ (the corresponding arguments are presented in section 4.3).
\begin{E}\label{notequal}
{\rm Assume that a closed subset $S\subset\partial\Di$ contains $-1$, $1$ and a sequence $\{e^{it_{k}}\}_{k\in\N}$, $t_{k}\in (0,\pi/2)$,  converging to $1$. Let $\{\alpha_{k}\}_{k\in\N}$ be a sequence of positive numbers satisfying the condition 
$\sum_{k=1}^{n}\alpha_{k}=1$. By $\chi_{k}$ we denote the characteristic function of the arc $\gamma_{k}:=\{e^{it}\ :\ t_{k}\leq t\leq\pi\}$. Consider the function 
$$
u(z):=\sum_{k=1}^{n}\alpha_{k}\chi_{k}(z),\ \ \ z\in\partial\Di.
$$
Clearly, $u\in R_{S}$. Let $h$ be a holomorphic function on $\Di$ such that
${\rm Re}\ \! h|_{\partial\Di}=u$. Then $e^{h}\in H^{\infty}\setminus A_{S}$.
However, for any $f\in A_{0}$ such that $f(1)=0$ we have $fe^{h}\in A_{S}$.
}
\end{E}
\begin{R}\label{lind}
{\rm
It seems to be natural that a function $f\in H^{\infty}$ has the weakest possible discontinuities on $\partial\Di$ if $f|_{\partial\Di}\in R_{S}$. However, 
from the classical Lindel\"{o}f theorem [L] it follows that any such $f$ in fact belongs to $A_{0}$. Moreover, the same conclusion is obtained even from the fact that ${\rm Re} f|_{\partial\Di}\in R_{S}$ for $f\in H^{\infty}$. In particular, if $f$ is holomorphic on $\Di$ and ${\rm Re} f|_{\partial\Di}$ is correctly defined and belongs to $R_{S}\setminus C(\partial\Di)$, then $f\not\in H^{\infty}$.\footnote{In this case $f|_{\partial\Di}\in BMO(\partial\Di)$ with $||f||_{BMO(\partial\Di)}\leq c||{\rm Re} f|_{\partial\Di}||_{L^{\infty}(\partial\Di)}$ for some absolute constant $c>0$.} Nevertheless, $e^{f}\in H^{\infty}$, which partly explains the choice of the object of our research.}
\end{R}

Let ${\cal M}(A_{S})$ be the maximal ideal space of $A_{S}$. Since evaluation
functionals $z(f):=f(z)$, $z\in\Di$, $f\in A_{S}$, belong to ${\cal M}(A_{S})$ and  $A_{S}$ separates points on $\Di$, there is a continuous embedding $i_{S}:\Di\to {\cal M}(A_{S})$. In the sequel we identify $\Di$ with $i_{S}(\Di)$. Then the following corona theorem is true.
\begin{Th}\label{te4}
$\Di$ is dense in ${\cal M}(A_{S})$.
\end{Th}
\begin{R}\label{r1}
{\rm Let us recall that the corona theorem is equivalent to the following statement, see, e.g., [G, Chapter V]:}

For any collection of functions $f_{1},\dots, f_{n}\in A_{S}$ satisfying
the corona condition
\begin{equation}\label{cc}
\max_{1\leq j\leq n}|f_{j}(z)|\geq\delta>0,\ \ \ z\in\Di,
\end{equation}
there are functions $g_{1},\dots, g_{n}\in A_{S}$ such that
\begin{equation}\label{bezout}
f_{1}g_{1}+\dots +f_{n}g_{n}=1.
\end{equation}
\end{R}

Finally we formulate some results on the structure of ${\cal M}(A_{S})$.
Since $A_{0}\hookrightarrow A_{S}$, there is a continuous
surjection of the maximal ideal spaces 
$$
a_{S}:{\cal M}(A_{S})\to {\cal M}(A_{0})\cong\overline{\Di}.
$$
Recall that the {\em \v{S}ilov boundary} of $A_{S}$ is the smallest compact subset $K\subset {\cal M}(A_{S})$ such that for each $f\in A_{S}$
$$
\sup_{z\in {\cal M}(A_{S})}|f(z)|=\sup_{\xi\in K}|f(\xi)|.
$$
Here we assume that every $f\in A_{S}$ is also defined on ${\cal M}(A_{S})$ where its extension to ${\cal M}(A_{S})\setminus\Di$ is given by the {\em Gelfand transform}: $f(\xi):=\xi(f)$, $\xi\in {\cal M}(A_{S})$.
\begin{Th}\label{te5}
\begin{itemize}
\item[(1)] $a_{S}:{\cal M}(A_{S})\setminus a_{S}^{-1}(S)\to\overline{\Di}\setminus S$ is a homeomorphism.
\item[(2)] The \v{S}ilov boundary $K_{S}$ of $A_{S}$ is naturally homeomorphic to ${\cal M}(AP(S))$. Under the identification of $K_{S}$ and ${\cal M}(AP(S))$ one has $a_{S}|_{K_{S}}=r_{S}\circ c_{S}$.
\item[(3)]
For each $z\in S$, the preimage $a_{S}^{-1}(z)$ is homeomorphic to the maximal ideal space of the algebra $AP_{{\cal O}}(\Sigma)$ of uniformly continuous almost periodic functions on the strip $\Sigma:=\{z\in\Co\ :\ {\rm Im}\ \!z\in [0,\pi]\}$ holomorphic at interior points of $\Sigma$.
\end{itemize}
\end{Th}

In the next section we describe the topological structure of the maximal ideal space ${\cal M}(AP_{{\cal O}}(\Sigma))$ of $AP_{{\cal O}}(\Sigma)$. We show that this space is equipped with a natural ``complex  structure``. Similarly each fibre $a_{S}^{-1}(z)$, $z\in S$, has a natural complex structure so that the homeomorphisms in Theorem \ref{te5}$\ \!$(3) are analytic with respect to these structures.

In a forthcoming paper we present similar results for algebras of bounded holomorphic functions on open polydisks generated by almost periodic functions.

\sect{\hspace*{-1em}. Maximal Ideal Space of the Algebra $AP_{{\cal O}}(\Sigma)$}
{\bf 2.1.} The construction presented below is rather general and can be defined for Galois coverings of complex manifolds with boundaries (cf. [Br]). However, we restrict ourselves to the case of coverings of annuli related to the subject of our paper.

Consider the annulus $R:=\{z\in\Co\ :\ e^{-2\pi^{2}}\leq |z|\leq 1\}$.
Its universal covering can be identified with $\Sigma$ so that
$e:\Sigma\to R$, $e(z):=e^{2\pi iz}$, $z\in\Sigma$, is the covering map.
We can also regard $\Sigma$ as a principal bundle on $R$ with fibre $\Z$
(see, e.g., [H] for the corresponding topological definitions). To specify,
consider a cover of $R$ by relatively open simply connected sets $U_{1}$ and $U_{2}$. Then $e^{-1}(U_{k})$ can be identified with $U_{k}\times\Z$, $k=1,2$.
Also, there is a continuous map $c_{12}:U_{1}\cap U_{2}\to\Z$ such that
$R$ is isomorphic (in the category of complex manifolds with boundaries) to the quotient space of $(U_{1}\times\Z)\sqcup (U_{2}\times\Z)$ under the equivalence relation: 
\begin{itemize}
\item[]
$U_{1}\times\Z\ni (z,n)\sim (z,n+c_{12}(z))\in U_{2}\times\Z$
for all $z\in U_{1}\cap U_{2}$ and $n\in\Z$. 
\end{itemize}
Let $b\Z$ be the Bohr compactification of $\Z$. Then the action 
of $\Z$ on itself by translations can be extended naturally to the action on $b\Z$: $\xi\mapsto \xi+n$, $\xi\in b\Z$, $n\in\Z$. By
$E(R,b\Z)$ we denote the principal bundle on $R$ with fibre $b\Z$ defined as the quotient of $(U_{1}\times b\Z)\sqcup (U_{2}\times b\Z)$ under the equivalence relation:
\begin{itemize}
\item[]
$U_{1}\times b\Z\ni (z,\xi)\sim (z,\xi+c_{12}(z))\in U_{2}\times b\Z$ for all $z\in U_{1}\cap U_{2}$ and $\xi\in b\Z$. 
\end{itemize}
Clearly $E(R,b\Z)$ is a compact Hausdorff space in the quotient topology induced by that of $(U_{1}\times b\Z)\sqcup (U_{2}\times b\Z)$. Also, the
projection $p: E(R,b\Z)\to R$ is defined by the natural projections
$U_{k}\times\Z\to U_{k}$, $k=1,2$, onto the first coordinate.

Next,
the natural local injections $U_{k}\times\Z\hookrightarrow U_{k}\times b\Z$, $k=1,2$, determine an injection $i_{0}:\Sigma\hookrightarrow E(R,b\Z)$ such 
that $p\circ i_{0}=e$. Moreover, $i_{0}(\Sigma)$ is dense in $E(R,b\Z)$, because $\Z$ is dense in $b\Z$ (in the topology of $b\Z$). Similarly, one  determines an injection $i_{\xi}:\Sigma\hookrightarrow E(R,b\Z)$, $\xi\in b\Z$, by the formula $i_{\xi}((z,n)):=(z,\xi+n)$,
$z\in U_{k}$, $n\in\Z$, $k=1,2$. Since, by definition, $\xi+\Z$ is dense in $b\Z$, the image $i_{\xi}(\Sigma)$ is dense in $E(R,b\Z)$ for any $\xi$.
Moreover, $E(R,b\Z)=\sqcup_{\xi}i_{\xi}(\Sigma)$ where
the union is taken over all $\xi$ whose images in the quotient group $b\Z/\Z$ are mutually distinct. Observe also that every $i_{\xi}$ is a continuous map and locally is an embedding.
\begin{D}\label{holom1}
A continuous function $f$ on $E(R,b\Z)$ is said to be holomorphic if  every function $f\circ i_{\xi}$ is holomorphic in interior points of $\Sigma$.
\end{D}

By ${\cal O}(E(R,b\Z))\subset C(E(R,b\Z))$ we denote the Banach algebra of holomorphic functions on $E(R,b\Z)$.
\begin{R}\label{holom2}
{\rm 
Using a normal family argument and the fact that $i_{\xi}(\Sigma)$ is dense in $E(R,b\Z)$ one can easily show that $f\in C(E(R,b\Z))$ is holomorphic if and only if there is $\xi\in b\Z$ such that $f\circ i_{\xi}$ is holomorphic in interior points of $\Sigma$.}
\end{R}
{\bf 2.2.} Let us recall that $f\in AP_{\cal O}(\Sigma)$ if $f$ is uniformly continuous with respect to the Euclidean metric on $\Sigma$, holomorphic in interior points of $\Sigma$ and its restriction to each straight line parallel to the $x$-axis is almost periodic. In the following lemma we identify $\Sigma$ with $i_{0}(\Sigma)\subset E(R,b\Z)$.
\begin{Lm}\label{holom3}
Every $f\in AP_{\cal O}(\Sigma)$ admits a continuous extension to a holomorphic function $\widehat f$ on $E(R,b\Z)$. Moreover, the correspondence 
$\widehat{}\ \!: AP_{\cal O}(\Sigma)\to {\cal O}(E(R,b\Z))$, $f\mapsto\widehat f$, determines an isomorphism of Banach algebras.
\end{Lm}
{\bf Proof.} Using the Pontryagin duality [P] one can easily show that the closure of $\Z$ in $b\Re$, the Bohr compactification of $\Re$, is isomorphic to $b\Z$. In particular, restrictions to $\Z$ of almost periodic functions on $\Re$ are almost periodic functions on $\Z$ and the algebra generated by extensions of such functions to $b\Z$ separates points on $b\Z$.
Let us identify naturally every $p^{-1}(U_{k})$ with $U_{k}\times b\Z$, and
every $e^{-1}(U_{k})$ with $U_{k}\times\Z\ (\subset U_{k}\times b\Z)$, $k=1,2$, and regard $f|_{e^{-1}(U_{k})}$, $f\in AP_{\cal O}(\Sigma)$, as a bounded function $f_{k}\in C(U_{k}\times\Z)$. Then, $f_{k}$ is
\begin{itemize}
\item[(a)]
holomorphic in interior points of $U_{k}\times\Z$, 
\item[(b)]
uniformly continuous on $U_{k}\times\Z$ with respect to the semi-metric $r(v_{1},v_{2}):=|z_{1}-z_{2}|$ on $U_{k}\times\Z$ where $v_{1}=(z_{1},n_{1})$, $v_{2}=(z_{2},n_{2})\in U_{k}\times\Z$ and 
$|\cdot|$ is the Euclidean norm on $\Co$, 
\item[(c)]
$f_{k}|_{\{z\}\times\Z}$ is almost periodic on $\Z$ for every $z\in U_{k}$.
\end{itemize}

To prove the lemma it suffices to show that 
\begin{itemize}
\item[(1)]
{\em
There is a continuous function $\widehat f_{k}$ on $U_{k}\times b\Z$  such
that $\widehat f_{k}|_{U_{k}\times\Z}=f_{k}$, for every $\xi\in b\Z$, the function $\widehat f_{k}|_{U_{k}\times\{\xi\}}$ is holomorphic in interior points of $U_{k}$ and $\sup_{U_{k}\times b\Z}|\widehat f_{k}|=\sup_{U_{k}\times\Z}|f_{k}|$.}
\item[(2)]
{\em  If $f\in {\cal O}(E(R,b\Z))$, then $f|_{\Sigma}\in AP_{\cal O}(\Sigma)$.}
\end{itemize}

(1) Since for every $z\in U_{k}$ the function $f_{kz}(n):=f_{k}(z,n)$ is almost periodic on $\Z$, there is a continuous function $\widehat f_{kz}$ on $b\Z$
which extends $f_{kz}$. We set 
$\widehat f_{k}(z,\xi):=\widehat f_{kz}(\xi)$, $\xi\in b\Z$ and prove that
$\widehat f_{k}$ is continuous. In fact, take a point 
$w=(z,\xi)\in U_{k}\times b\Z$ and a number $\epsilon>0$.
By the uniform continuity of $f_{k}$, there is $\delta>0$ such that for any pair of points
$v_{1}=(z_{1},n)$ and $v_{2}=(z_{2},n)$ from $U_{k}\times\Z$ with
$|z_{1}-z_{2}|<\delta$ one has $|f(z_{1},n)-f(z_{2},n)|<\epsilon/3$. Define a
neighbourhood $U_{z}$ of $z\in U_{k}$ by $U_{z}:=\{z'\in U_{k}\ : |z-z'|<\delta\}$.
Further, by the definition of $\widehat f_{kz}$, there is a
neighbourhood $U_{\xi}\subset b\Z$ of $\xi$ such that for any $\eta\in
U_{\xi}$ we have $|\widehat f_{kz}(\eta)-\widehat f_{kz}(\xi)|<\epsilon/3$.
Consider $U_{w}:=U_{z}\times U_{\xi}$. Then $U_{w}$ is
an open neighbourhood of $w\in U_{k}\times b\Z$. Note that
$f_{kz}-f_{kz'}$ is an almost periodic function on $\Z$ and 
for any $z'\in U_{z}$ its supremum norm is $<\epsilon/2$. This implies that
$|\widehat f_{kz}(\eta)-\widehat f_{kz'}(\eta)|<\epsilon/2$ for each
$\eta\in b\Z$. In particular, for any $(x,\eta)\in U_{w}$ we have
$$
|\widehat f_{k}(x,\eta)-\widehat f_{k}(z,\xi)|\leq |\widehat f_{kx}(\eta)-\widehat
f_{kz}(\eta)|+|\widehat f_{kz}(\eta)-\widehat f_{kz}(\xi)|<\epsilon \ .
$$
This shows that $\widehat f_{k}$ is continuous at every
$w\in U_{k}\times b\Z$.

Now, we show that $\widehat f|_{U_{k}\times\{\xi\}}$ is holomorphic in interior points of $U_{k}$ for every $\xi\in b\Z$.

Since $\widehat f_{k}$ is uniformly continuous on the compact set
$U_{k}\times b\Z$, for any $\epsilon>0$ there is 
$n_{\epsilon}\in\Z$ such that 
$\sup_{z\in U_{k}}|\widehat f_{k}(z,\xi)-f_{k}(z,n_{\epsilon})|<\epsilon$. 
In particular,
$\widehat f_{k}(\cdot,\xi)$ is the limit in $C(U_{k})$ of the sequence 
$\{f_{k}(\cdot,n_{1/l})\}_{l\geq 1}$
of bounded continuous functions holomorphic on the interior of $U_{k}$. Thus $\widehat f_{k}|_{U_{k}\times\{\xi\}}$ is also holomorphic on the interior of $U_{k}$.

Note that the equality $\sup_{U_{k}\times b\Z}|\widehat f_{k}|=\sup_{U_{k}\times\Z}|f_{k}|$ follows directly from the definition of $\widehat f_{k}$. This completes the proof of (1).

(2) Suppose that $f\in {\cal O}(E(R,b\Z))$. We must show that $f|_{\Sigma}\in AP_{\cal O}(\Sigma)$. To this end it suffices to show that for every line
$L:=\{z\in\Co\ :\ {\rm Im}\ \!z=t\in [0,\pi]\}$, the function $f|_{L}$ is almost periodic. (The uniform continuity of $f|_{\Sigma}$ with respect to the Euclidean metric on $\Co$ follows easily from the uniform continuity of $f$ on $E(R,b\Z)$.) By definition the image $S:=e(L)\subset R$ is a circle and  $e|_{L}:L\to S$ is the universal covering.  Consider the compact set $p^{-1}(S)\subset E(R,b\Z)$. Since the function $f|_{p^{-1}(S)}$ is continuous and $b\Z\subset b\Re$, given $\epsilon>0$ there are a finite open cover $(V_{n})_{1\leq n\leq m}$ of $S$ by sets homeomorphic to open intervals in $\Re$ and continuous almost periodic functions $f_{n}$ on $L$, $1\leq n\leq m$, such that
\begin{equation}\label{approxim1}
\sup_{z\in e^{-1}(V_{n}),1\leq n\leq m}|f_{n}(z)-f(z)|<\epsilon .
\end{equation}
Let $\{\rho_{n}\}_{1\leq n\leq m}$ be a continuous partition of unity subordinate to $(V_{n})_{1\leq n\leq m}$. We pull back it to $L$ by $e$ and by $\widetilde\rho_{n}\ (:=e^{*}\rho_{n})$, $1\leq n\leq m$, denote the obtained functions. Since each $\widetilde\rho_{n}$ is periodic on $L$, it is almost periodic, as well. Let us define the function $f_{\epsilon}$ on $L$ by the formula
$$
f_{\epsilon}(z):=\sum_{n=1}^{m}\widetilde\rho_{n}(z)f_{n}(z).
$$
Then clearly $f_{\epsilon}$ is a continuous almost periodic function on $L$
and by (\ref{approxim1})
$$
\sup_{z\in L}|f_{\epsilon}(z)-f(z)|<\epsilon.
$$
This shows that $f$ admits uniform approximation on $L$ by continuous almost periodic functions and therefore $f|_{L}$ is almost periodic, as well.

The lemma has been proved.\ \ \ \ \ $\Box$\\
\\
{\bf 2.3.} In this part we prove the corona theorem for the algebra $AP_{\cal O}(\Sigma)$.
Recall that ${\cal M}(AP_{\cal O}(\Sigma))$ stands for its maximal ideal space. It is well known that every $f\in AP_{\cal O}(\Sigma)$ can be approximated uniformly on $\Sigma$ by polynomials in $e^{i\lambda z}$, $\lambda\in\Re$, see, e.g., [JT]. Then using the inverse limit construction for maximal ideal spaces of uniform algebras, see [R], one obtains that the base of topology of ${\cal M}(AP_{\cal O}(\Sigma))$ is generated by functions $e^{i\lambda z}$. Namely, the base of topology on ${\cal M}(AP_{\cal O}(\Sigma))$ consists of open sets of the form
$$
U(\lambda_{1},\dots,\lambda_{l},\xi,\epsilon):=\{\eta\in {\cal M}(AP_{\cal O}(\Sigma))\ :\ \max_{1\leq k\leq l}|e_{\lambda_{k}}(\eta)-e_{\lambda_{k}}(\xi)|<\epsilon\}
$$
where $e_{\lambda}$ is the extension of $e^{i\lambda z}$ to ${\cal M}(AP_{\cal O}(\Sigma))$ by means of the Gelfand transform. 
\begin{Th}\label{coron1}
$\Sigma$ is dense in ${\cal M}(AP_{\cal O}(\Sigma))$ in the Gelfand topology.
\end{Th}
{\bf Proof.}
Assume that the corona theorem is not true for $AP_{\cal O}(\Sigma)$, that is, $\Sigma$ is not dense in ${\cal M}(AP_{\cal O}(\Sigma))$. Then there exist $\xi\in {\cal M}(AP_{\cal O}(\Sigma))$ and its neighbourhood $U(\lambda_{1},\dots,\lambda_{l},\xi,\epsilon)$ such that
$$
U(\lambda_{1},\dots,\lambda_{l},\xi,\epsilon)\cap cl(\Sigma)=\emptyset;
$$
here $cl(\Sigma)$ is the closure of $\Sigma$ in ${\cal M}(AP_{\cal O}(\Sigma))$.
Denoting $c_{k}:=e_{\lambda_{k}}(\xi)$, $1\leq k\leq l$, we have
\begin{equation}\label{strip1}
\max_{1\leq k\leq l}|e^{i\lambda_{k}z}-c_{k}|\geq\epsilon>0\ \ \ {\rm for\ all}\ \ \ z\in\Sigma.
\end{equation}
Clearly, every function $e^{i\lambda_{k}z}-c_{k}$, $1\leq k\leq l$, has at least one zero in $\Sigma$. (For otherwise, if, say,  $e^{i\lambda_{k}z}-c_{k}$ has no zeros on $\Sigma$, then the function $g_{k}(z):=\frac{1}{e^{i\lambda_{k}z}-c_{k}}$, $z\in\Sigma$, obviously belongs to $AP_{\cal O}(\Sigma)$ and $g_{k}(z)(e^{i\lambda_{k}z}-c_{k})=1$ for all $z\in\Sigma$, a contradiction to our assumption.) In particular, since solutions of the equation
$$
e^{i\lambda_{k}z}=c_{k},\ \ \ \lambda_{k}\neq 0,
$$
are given by
$$
z=-\frac{i\ln|c_{k}|}{\lambda_{k}}+\frac{{\rm Arg}\ \! c_{k}+2\pi s}{\lambda_{k}},\ \ \ s\in\Z,
$$ 
they all belong to $\Sigma$. Further, without loss of generality we may assume that all $\lambda_{k}>0$. Indeed, if some $\lambda_{k}<0$, we can replace the function $e^{i\lambda_{k}z}-c_{k}$ by $e^{-i\lambda_{k}z}-\frac{1}{c_{k}}$ (observe that $c_{k}\neq 0$ by the above argument) so that the new family of functions also satisfies  (\ref{strip1}) (possibly with a different $\epsilon$) and extensions of these functions to ${\cal M}(AP_{\cal O}(\Sigma))$ vanish at $\xi$. Since all these functions have zeros in $\Sigma$ and satisfy (\ref{strip1}) there, we have
$$
\max_{1\leq k\leq l}|e^{i\lambda_{k}z}-c_{k}|\geq\widetilde\epsilon>0,\ \ \ {\rm for\ all}\ \ \ z\in\H_{+}
$$
where $\H_{+}\subset\Co$ is the open upper half-plane, and all $e^{i\lambda_{k}z}-c_{k}$, $\lambda_{k}\in\Re_{+}$, are almost periodic on $\overline\H_{+}$. From the last inequality by the B\"{o}tcher corona theorem [B\"{o}] we obtain that there exist holomorphic almost periodic functions $g_{1},\dots, g_{l}$ on $\overline \H_{+}$ such that
$$
\sum_{k=1}^{l}g_{k}(z)(e^{i\lambda_{k}z}-c_{k})=1\ \ \ {\rm for\ all}\ \ \
z\in\overline\H_{+}.
$$
Thus taking the restrictions of these functions to $\Sigma$ we obtain a contradiction to our assumption.

This completes the proof of the corona theorem for $AP_{\cal O}(\Sigma)$.\ \ \ \ \ $\Box$
\begin{C}\label{corona2}
$E(R,b\Z)$ is homeomorphic to ${\cal M}(AP_{\cal O}(\Sigma))$.
\end{C}
{\bf Proof.} By Lemma \ref{holom3} there exists a continuous embedding $i$ of $E(R,b\Z)$ into ${\cal M}(AP_{\cal O}(\Sigma))$. Since $\Sigma$ is dense in $E(R,b\Z)$ and $i(\Sigma)$ is dense in ${\cal M}(AP_{\cal O}(\Sigma))$ by the corona theorem, $i$ is a homeomorphism.\ \ \ \ \ $\Box$
\begin{R}\label{bound1}
{\rm (1) According to our construction, the closure of $\Re$ in $E(R,b\Z)$ coincides with $b\Re$. In fact, this closure is $E(\partial\Di,b\Z)$, the principal bundle on $\partial\Di$ with fibre $b\Z$ obtained as the restriction of $E(R,b\Z)$ to $\partial\Di$. Since $R$ is homotopically equivalent to $\partial\Di$, by the covering homotopy theorem $E(R,b\Z)$ is homotopically equivalent to $b\Re=E(\partial\Di,b\Z)$.\\
(2) Observe also that the covering dimension of $b\Re$ is $\infty$, because $b\Re$ is the inverse limit of real tori whose dimensions go to $\infty$. Therefore the covering dimension of ${\cal M}(AP_{\cal O}(\Sigma))=E(R,b\Z)$ is also $\infty$.\\
(3) Finally, it is easy to show that the \v{S}ilov boundary of $AP_{\cal O}(\Sigma)$ is $E(R,b\Z)|_{\partial R}$, the restriction of $E(R,b\Z)$ to the boundary $\partial R$ of $R$, and is homeomorphic to $b\Re\sqcup b\Re$.}
\end{R}  
\sect{\hspace*{-1em}. Proofs of Theorems \ref{te1}, \ref{te2} and Corollary \ref{cor1}}
{\bf 3.1. Proof of Theorem \ref{te1}.} 
Let $f\in AP(S)$.
Since $\partial\Di$ is a compact set, given $\epsilon>0$ there are finitely many points $z_{l}:=e^{it_{l}}\in \partial\Di$, numbers $s_{l}\in (0,\pi)$ and almost periodic functions $f_{k}^{l}:\gamma_{t_l^{k}}(s_l) \to \Co$, $k\in\{-1,1\}$, $1 \leq l \leq n$, such that 
\begin{equation}\label{e2}
\begin{array}{c}
\displaystyle
\bigcup_{l=1}^{n}(\gamma_{t_{l}^{-1}}(s_{l})\cup \gamma_{t_{l}^{1}}(s_{l})) =\partial\Di\setminus\{z_{1},\dots, z_{n}\} 
\ \ \ {\rm and\ for\ all}\ \ \ 1\leq l\leq n
\\
\\
\displaystyle
{\rm ess}\!\!\!\!\!\sup_{z \in \gamma_{t_{l}^{1}}(s_l)}|f(z)-f^l_{1}(z)|<\frac{\epsilon}{2}\ ,\ \ \ \ 
\ \ \ {\rm ess}\!\!\!\!\!\sup_{z \in \gamma_{t_{l}^{-1}}(s_l)}|f(z)-f^{l}_{-1}(z)|<\frac{\epsilon}{2}\ .
\end{array}
\end{equation}
We set $U_{l}:=\gamma_{t_{l}^{1}}(s_{l})\cup \gamma_{t_{l}^{-1}}(s_{l})\cup\{z_{l}\}$.
Then $U=(U_{l})_{l=1}^{n}$ is a finite open cover of $\partial\Di$. Let $\{\rho_{l}\}_{l=1}^{n}$ be a continuous partition of unity subordinate to $U$ such that ${\rm supp}\ \!\rho_{l}\subset\subset U_{l}$ and $\rho_{l}(z_{l})=1$, $1\leq l\leq n$.
Consider the functions $f_{l}$ on $\partial\Di\setminus\{z_{l}\}$ defined by the 
formulas
$$
f_{l}(z):=\left\{
\begin{array}{ccc}
\rho_{l}(z)f_{-1}^{l}(z),&{\rm if}&z\in \gamma_{t_{l}^{-1}}(s_{l}),\\
\\
\rho_{l}(z)f_{1}^{l}(z),&{\rm if}&z\in \gamma_{t_{l}^{1}}(s_{l}).
\end{array}
\right.
$$
Since $f_{l}$ is continuous outside $z_{l}$ and coincides with $f_{-1}^{l}$ and $f_{1}^{l}$ in a neighbourhood of $z_{l}$, it belongs to $AP(S)$. Moreover,
\begin{equation}\label{e21}
||f-\sum_{l=1}^{n}f_{l}||_{L^{\infty}(\partial\Di)}<\frac{\epsilon}{2}.
\end{equation}
Thus in order to prove the theorem it suffices to approximate every  $f_{l}$ by polynomials in $g_{t_{l}}$ and
$e^{i\lambda\log_{t_{l}^{k}}}$, $k\in\{-1,1\}$, $\lambda\in\Re$. 

Suppose first that $z_{l}\not\in S$. Since $f$ is continuous outside the compact set $S$, we can choose the above cover $U$ and the family of functions $\{f_{k}^{m}\}_{1\leq m\leq n}$, $k\in\{-1,1\}$, such that in the $U_{l}$ the functions $f_{k}^{l}$,
$k\in\{-1,1\}$, have the same limit at $z_{l}$. This implies the continuity of $f_{l}$ on $\partial\Di$. Next,
consider the uniform algebra $\Co(g_{t_{l}})$ over $\Co$ generated by the function $g_{t_{l}}$. Since by our definition $g_{t_{l}}$ separates points on $\partial\Di\setminus\{z_{l}\}$, the maximal ideal space of $\Co(g_{t_{l}})$ is homeomorphic to the closed interval $(\partial\Di\setminus\{z_{l}\})\cup\{z_{{l}^{k}}\}$, $k\in\{-1,1\}$, with endpoints  
$z_{{l}^{1}}$ and $z_{{l}^{-1}}$ identified with the counterclockwise and clockwise orientations at $z_{l}$.
Clearly every continuous function on $\partial\Di$ is extended to the maximal ideal space of $\Co(g_{t_{l}})$ as a continuous function having the same values at $z_{{l}^{1}}$ and $z_{{l}^{-1}}$. Thus by the Stone-Weierstrass theorem $f_{l}$ can be uniformly approximated on $\partial\Di\setminus\{z_{l}\}$ by complex polynomials in $g_{t_{l}}$.

Now, suppose that $z_{l}\in S$. Choose some $s\in (0,s_{l})$.
By $AP_{\{z_{l}\}}(s)$ we denote the uniform algebra of complex continuous functions on $\partial\Di\setminus\{z_{l}\}$ almost periodic on the open arcs $\gamma_{t_{l}^{k}}(s)$, $k\in\{-1,1\}$. (Since $s\in (0,\pi)$, the closures of these arcs are disjoint.) By ${\cal M}_{\{z_{l}\}}(s)$ we denote the maximal ideal space of $AP_{\{z_{l}\}}(s)$.  Then $\partial\Di\setminus\{z_{l}\}$ is dense in ${\cal M}_{\{z_{l}\}}(s)$ (in the Gelfand topology). Note that the space ${\cal M}_{\{z_{l}\}}(s)$ is constructed as follows.

Consider the Bohr compactification $b\Re$ of $\Re$. We identify the negative ray $\Re_{-}$ in $\Re\subset b\Re$ with $\gamma_{t_{{l}^{1}}}(s)\subset \partial\Di$ by means of the map $t\mapsto e^{i(t_{l}+se^{t})}$, $t\in\Re_{-}$. Similarly, consider another copy of  $b\Re$ and identify $\Re_{-} (\subset\Re)$ in this copy with $\gamma_{t_{{l}^{-1}}}(s)\subset \partial\Di$ by means of the map $t\mapsto e^{i(t_{l}-se^{t})}$, $t\in\Re_{-}$. On the identified sets we introduce the topology induced from 
$b\Re$ and on $\partial\Di\setminus(\gamma_{t_{{l}^{-1}}}(s)\cup \gamma_{t_{{l}^{1}}}(s))$ the  topology is induced from $\partial\Di$. Then the quotient space of $b\Re\sqcup b\Re\sqcup\partial\Di$ under these identifications coincides with ${\cal M}_{\{z_{l}\}}(s)$.

Next, recall that by definition the algebra $AP(\{z_{l}\})$ is the uniform closure in $C(\partial\Di\setminus\{z_{l}\})$ of the algebra generated by the algebras $AP_{\{z_{l}\}}(s)$, $s\in (0,s_{l})$. By ${\cal M}(AP(\{z_{l}\}))$ we denote the maximal ideal space of $AP(\{z_{l}\})$. Since for any $s''<s'$ we have inclusions
$i_{s''s'}:AP_{\{z_{l}\}}(s')\hookrightarrow AP_{\{z_{l}\}}(s'')$, the space ${\cal M}(AP(\{z_{l}\}))$ is the inverse limit of the spaces ${\cal M}_{\{z_{l}\}}(s)$ (here the maps
$p_{s''s'}:{\cal M}_{\{z_{l}\}}(s'')\to {\cal M}_{\{z_{l}\}}(s')$ in the definition of this limit are defined as the dual maps to $i_{s''s'}$). Also, $\partial\Di\setminus\{z_{l}\}$ is dense in ${\cal M}(AP(\{z_{l}\}))$ in the Gelfand topology. Since by the definition the functions $f_{l}$, $e^{i\lambda\log_{t_{{l}^{k}}}}$ and $g_{t_{l}}$ admit the continuous extensions (denoted by the same symbols) to ${\cal M}(AP(\{z_{l}\}))$, it suffices to show that the extended functions $e^{i\lambda\log_{t_{{l}^{k}}}}$, $g_{t_{l}}$ separate points on ${\cal M}(AP(\{z_{l}\}))$. Then we will apply the Stone-Weierstrass theorem to get a complex polynomial $p_{l}$ in variables $e^{i\lambda\log_{{t_{l}^{k}}}}$ and $g_{t_{l}}$ which uniformly approximates $f_{l}$ on ${\cal M}(AP(\{z_{l}\}))$ with an error $<\epsilon/2n$. Therefore, $\sum_{l=1}^{n}p_{l}$ will approximate $f$ in
$L^{\infty}(\partial\Di)$ with an error $<\epsilon$. 

So let us show that the functions $e^{i\lambda\log_{{t_{l}^{k}}}}$, $g_{t_{l}}$ separate points on ${\cal M}(AP(\{z_{l}\}))$. 
By $p_{s}:{\cal M}(AP(\{z_{l}\}))\to {\cal M}_{\{z_{l}\}}(s)$ we denote the continuous surjection determined by the inverse limit construction.
First, we consider distinct points $x,y\in {\cal M}(AP(\{z_{l}\}))$ for which there is $s\in (0,s_{l})$ such that  $p_{s}(x)$ and $p_{s}(y)$ are distinct and belong to one of the Bohr compactifications of $\Re$ in ${\cal M}_{\{z_{l}\}}(s)$, say, e.g., to the compactification obtained by gluing $\Re_{-}$ with $\gamma_{t_{{l}^{1}}}(s)$. Since in this case the functions $e^{i\lambda\log_{t_{{l}^{1}}}}$ extended to $b\Re$ are identified with the extensions to $b\Re$ of functions $c_{\lambda s}e^{i\lambda t}$ on $\Re$, 
$c_{\lambda s}:=e^{i\lambda\ln s}$, by the classical Bohr theorem there is $\lambda_{0}\in\Re$ such that the extension of $e^{i\lambda_{0}\log_{t_{{l}^{1}}}}$ to $b\Re$ separates $p_{s}(x)$ and $p_{s}(y)$. So, the extension of $e^{i\lambda_{0}\log_{t_{{l}^{1}}}}$ to
${\cal M}(AP(\{z_{l}\}))$ separates $x$ and $y$.

Suppose now that $x$ and $y$ are such that $p_{s}(x)$ and $p_{s}(y)$ belong to different Bohr compactifications of $\Re$ for all $s$. This implies that 
$x$ and $y$ are limit points of the sets $\gamma_{t_{l}^{k(x)}}(s)$ and $\gamma_{t_{l}^{k(y)}}(s)$ for some $s\in (0,s_{l})$ with $k(x)\neq k(y)$ and $k(x),k(y)\in\{-1,1\}$.
Then the function $g_{t_{l}}$ by definition equals 1 at one of the points $x$, $y$ and 0 at the other one.

Finally, assume that $x\in {\cal M}(AP(\{z_{l}\}))\setminus\partial\Di$ and $y\in\partial\Di\setminus\{z_{l}\}$. Then $g_{t_{l}}(x)$ equals either $0$ or $1$ and $g_{t_{l}}(y)$ differs from these numbers because $g_{t_{l}}$ is decreasing on $\partial\Di\setminus\{z_{l}\}$.

Thus we have proved that the family of functions $e^{i\lambda\log_{{t_{l}^{k}}}}$, $g_{t_{l}}$ separate points on ${\cal M}(AP(\{z_{l}\}))$. This completes the proof of the theorem.\ \ \ \ \ $\Box$\\
\\
{\bf 3.2. Proof of Corollary \ref{cor1}.} According to Theorem \ref{te1} it suffices to prove that for functions $g_{t_{0}}$ and $e^{i\lambda\log_{t_{0}^{k}}}$, $\lambda\in\Re$, $z_{0}:=e^{it_{0}}\in\widetilde S:=\phi(S)$, $k\in\{-1,1\}$, the corresponding functions $\phi^{*}(g_{t_{0}})$ and 
$f:=\phi^{*}(e^{i\lambda\log_{t_{0}^{k}}})$ belong to $AP(S)$. Since $g_{t_{0}}\in R_{\{z_{0}\}}$, the statement is trivial for $\phi^{*}(g_{t_{0}})$.  
Without loss of generality we may assume that $\phi$ preserves the orientation on $\partial\Di$. Let $e^{i\tilde t_{0}}:=\phi^{-1}(z_{0})$. Then we have
$$
\phi(e^{i(\tilde t_{0}+t)})=e^{i(t_{0}+\tilde\phi(t))},\ \ \ t\in [0,2\pi],
$$
where $\tilde\phi:[0,2\pi]\to [0,2\pi]$ is a $C^{1}$ diffeomorphism and
$\tilde\phi(0)=0$. We set $c:=\tilde \phi'(0)$. 
By the definition on the open arcs $\gamma_{\tilde t_{0}^{k}}(s)$ for
$t\in (-\infty,0)$ we obtain
$$
\widehat f(t):=f(e^{i(\tilde t_{0}+kse^{t})})=
e^{i\lambda\ln\tilde\phi(se^{t}))}=e^{i\lambda\ln(cse^{t}+o(se^{t}))}=e^{i(\lambda\ln(cs)+o(1))}e^{i\lambda t},\ \ \ {\rm as}\ \ \ s\to 0.
$$
Since $f$ is continuous outside $e^{i\tilde t_{0}}$, the latter implies that
$f\in AP(S)$.\ \ \ \ \ $\Box$
\\
\\
{\bf 3.3. Proof of Theorem \ref{te2}.} We begin with the following 
\begin{Lm}\label{unif1}
The algebra $R_{S}$ is the uniform closure of the algebra generated by the algebras $R_{\{z\}}$ for all possible $z\in S$.
\end{Lm}
(For $S=\partial\Di$ a similar statement first was proved by Dieudonne [D]).\\
{\bf Proof.}
Consider a regulated function $f\in R_{S}$. Since $\partial\Di$ is a compact set, by the definition of $R_{S}$ for any $\epsilon>0$ there are finitely many points $z_{l}:=e^{it_{l}}\in\partial\Di$, numbers $s_{l}\in (0,\pi)$
and constant functions $f_{k}^{l}:\gamma_{t_l^{k}}(s_l) \to \Co$, $k\in\{-1,1\}$, $1 \leq l \leq n$, such that 
$$
\begin{array}{c}
\displaystyle
\bigcup_{l=1}^{n}(\gamma_{t_{l}^{-1}}(s_{l})\cup \gamma_{t_{l}^{1}}(s_{l}))=\partial\Di\setminus\{z_{1},\dots, z_{n}\} 
\ \ \ {\rm and\ for\ all}\ \ \ 1\leq l\leq n
\\
\\
\displaystyle
{\rm ess}\!\!\!\!\!\sup_{z \in \gamma_{t_{l}^{1}}(s_l)}|f(z)-f^l_{1}(z)|<\epsilon\ ,\ \ \ \ 
\ \ \ {\rm ess}\!\!\!\!\!\sup_{z \in \gamma_{t_{l}^{-1}}(s_l)}|f(z)-f^{l}_{-1}(z)|<\epsilon\ .
\end{array}
$$
We set $U_{l}:=\gamma_{t_{l}^{1}}(s_{l})\cup \gamma_{t_{l}^{-1}}(s_{l})\cup\{z_{l}\}$.
Then $U=(U_{l})_{l=1}^{n}$ is a finite open cover of $\partial\Di$. Let $\{\rho_{l}\}_{l=1}^{n}$ be a continuous partition of unity subordinate to $U$ such that ${\rm supp}\ \!\rho_{l}\subset\subset U_{l}$ and $\rho_{l}(z_{l})=1$, $1\leq l\leq n$.
Consider the functions $f_{l}$ on $\partial\Di\setminus\{z_{l}\}$ defined by the 
formulas
$$
f_{l}(z):=\left\{
\begin{array}{ccc}
\rho_{l}(z)f_{-1}^{l}(z),&{\rm if}&z\in \gamma_{t_{l}^{-1}}(s_{l}),\\
\\
\rho_{l}(z)f_{1}^{l}(z),&{\rm if}&z\in \gamma_{t_{l}^{1}}(s_{l}).
\end{array}
\right.
$$
If $z_{l}\in S$, then by definition $f_{l}\in R_{\{z_{l}\}}$.
If $z_{l}\not\in S$, then since $f$ is continuous outside the compact set $S$, we can choose the above cover $U$ and the family of functions $\{f_{k}^{m}\}_{1\leq m\leq n}$, $k\in\{-1,1\}$, such that in the $U_{l}$ the functions $f_{k}^{l}$,
$k\in\{-1,1\}$, have the same limit at $z_{l}$. This implies the continuity of $f_{l}$ on $\partial\Di$, i.e., $f_{l}\in R_{\{z_{l}\}}$, as well. Also, we have
$$
||f-\sum_{l=1}^{n}f_{l}||_{L^{\infty}(\partial\Di)}<\epsilon.
$$

This completes the proof of the lemma.\ \ \ \ \ $\Box$

Let $F_{1}\subset F_{2}$ be finite subsets of $S$. Then we have the natural injection $i_{F_{1}F_{2}}:R_{F_{1}}\hookrightarrow R_{F_{2}}$. Passing to the
map dual to $i_{F_{1}F_{2}}$ we obtain a continuous map of the corresponding
maximal ideal spaces: $p_{F_{1}F_{2}}:{\cal M}(R_{F_{2}})\to {\cal M}(R_{F_{1}})$. Now the family $\{({\cal M}(R_{F_{1}}),{\cal M}(R_{F_{2}}),
p_{F_{1}F_{2}})\}_{F_{1}\subset F_{2}\subset S}$ determines an inverse limit system whose limit according to Lemma \ref{unif1} coincides with ${\cal M}(R_{S})$. By $p_{F}: {\cal M}(R_{S})\to {\cal M}(R_{F})$, $F\subset S$ is finite, we denote the limit maps determined by this limit. Then we have
\begin{equation}\label{comp}
c_{F}\circ p_{F}=c_{S}.
\end{equation}

Let $F$ consist of $n$ points. Then $\partial\Di\setminus F$ is a disjoint union of open arcs $\gamma_{k}$, $1\leq k\leq n$.
Consider a real function $g_{F}$ on $\partial\Di$ having discontinuities of the first kind at points $F$ and continuous outside $F$ such that the  $g_{F}:\partial\Di\setminus F\to\Re$ is an injection and $g_{F}(\partial\Di\setminus F)$ is the union of open intervals whose closures are mutually disjoint. Then an argument similar to that used in the proof of Theorem \ref{te1} (see the case $z_{l}\notin S$ there) shows that 
$R_{F}=\Co(g_{F})$, the uniform algebra in $L^{\infty}(\partial\Di)$ generated by $g_{F}$. This implies that ${\cal M}(R_{F})$ is naturally homeomorphic to the disjoint union of closures $\overline \gamma_{k}$ of $\gamma_{k}$, $1\leq k\leq n$, and $c_{F}:{\cal M}(R_{F})\to\partial\Di$ maps identically every $\gamma_{k}$ in this union to $\gamma_{k}\subset\partial\Di$. Since $c_{F}$ is continuous, for every $z\in F$ the preimage $c_{F}^{-1}(z)$ consists of  two points $z_{+}$ and $z_{-}$ that can be naturally identified with counterclockwise and clockwise orientations of $\partial\Di$ at $z$.
Thus we obtain proofs of statement (1)-(3) of the theorem for a finite subset $F\subset S$.  To prove the general case we use (\ref{comp}).
 
Assume that for some $z\in\partial\Di$ the preimage $c_{S}^{-1}(z)$ contains at least three points $x_{1}$, $x_{2}$ and $x_{3}$. Then by the definition of the inverse limit, there is a finite subset $F\subset S$ such that
$p_{F}(x_{1})$, $p_{F}(x_{2})$ and $p_{F}(x_{3})$ are distinct points in
${\cal M}(R_{F})$. Since by (\ref{comp}) the images of these points under $c_{F}$ coincide with $z$, from the case established above for ${\cal M}(R_{F})$ we obtain that $c_{F}^{-1}(z)$ consists of at most two points, a contradiction. 
Thus $c_{S}^{-1}(z)$ consists of at most 2 points for every $z\in\partial\Di$.

Assume now that $z\in S$.  Let $F\subset S$ be a finite subset containing $z$.
Then the preimage $c_{F}^{-1}(z)$ consists of two points. Since by (\ref{comp})
$c_{S}^{-1}(z)=p_{F}^{-1}(c_{F}^{-1}(z))$, the preimage $c_{S}^{-1}(z)$ consists of two points, as well. (As before, they can be naturally identified with counterclockwise and clockwise orientations of $\partial\Di$ at $z$.)

This proves statement (1) of the theorem.

Next, for $z\not\in S$ we have $c_{S}^{-1}(z)$ is a single point. (For otherwise, for some finite $F$ we will have that $c_{F}^{-1}(z)$ consists of at least 2 points, a contradiction.) Since $c_{S}^{-1}(S)$ is compact, the latter implies that $c_{S}:{\cal M}(R_{S})\setminus c_{S}^{-1}(S)\to\partial\Di\setminus S$ is a homeomorphism. 

The proof of statement (2) is complete.

To prove (3) we will assume that $S$ is infinite (for finite $S$ the statement is already proved). Let $F\subset S$ be a finite subset consisting of $n$ points, $n\geq 2$. Let $\partial\Di\setminus F$ be the disjoint union of 
open arcs $\gamma_{k}$, $1\leq k\leq n$. By $\chi_{\gamma_{k}}$ we denote the characteristic function of $\gamma_{k}$. Then every $\chi_{\gamma_{k}}$ belongs to $R_{F}$. By the same symbols we denote continuous extensions of $\chi_{\gamma_{k}}$ to ${\cal M}(R_{F})$ by means of the Gelfand transform. Let us determine a continuous map $K_{F}:{\cal M}(R_{F})\to \Z_{2}(F):=\{0,1\}^{n}$ by the formula
$$
K_{F}(m):=(\chi_{\gamma_{1}}(m),\dots,\chi_{\gamma_{n}}(m)),\ \ \ m\in {\cal M}(R_{F}).
$$
Set ${\cal Z}(S):=\prod_{F\subset S}\Z_{2}(F)$ and determine the map ${\cal K}_{S}:{\cal M}(R_{S})\to {\cal Z}(S)$ by the formula
$$
{\cal K}_{S}(m):=\{K_{F}(p_{F}(m))\}_{F\subset S}.
$$
We equip ${\cal Z}(S)$ with the {\em Tychonoff topology}. Then ${\cal Z}(S)$ is a totally disconnected compact Hausdorff space and the map ${\cal K}_{S}$ is continuous. 

Show that ${\cal K}_{S}|_{c_{S}^{-1}(S)}:c_{S}^{-1}(S)\to {\cal Z}(S)$ is an injection. Indeed, for distinct points $x$, $y$ from $c_{S}^{-1}(S)$ there exists a finite subset $F\subset S$ consisting of at least two points such that $p_{F}(x)\neq p_{F}(y)$ and $p_{F}(x), p_{F}(y)\in c_{F}^{-1}(F)$.
Then by the definition of the map $K_{F}$ we have
$$
K_{F}(p_{F}(x))\neq K_{F}(p_{F}(y)).
$$
This implies that ${\cal K}_{S}(x)\neq {\cal K}_{S}(y)$. Since $c_{S}^{-1}(S)$ is a compact set, the injectivity implies that $c_{S}^{-1}(S)$ is homeomorphic to ${\cal K}_{S}(c_{S}^{-1}(S))$.
The latter space is totally disconnected as a compact subset of the totally disconnected space ${\cal Z}(S)$.

This completes the proof of (3).

(4) According to (\ref{e21}) the maximal ideal space ${\cal M}(AP(S))$ of the algebra $AP(S)$ is homeomorphic to the inverse limit of compact spaces ${\cal M}(AP(F))$ with $F\subset S$ finite.
Let $\widetilde p_{F_{1}F_{2}}:{\cal M}(AP(F_{2}))\to {\cal M}(AP(F_{1}))$, $F_{1}\subset F_{2}$, be continuous maps determining this limit, and $\widetilde p_{F}:{\cal M}(AP(S))\to {\cal M}(AP(F))$ be the corresponding limit maps. Since each $AP(F)$ is a self-adjoint algebra, by the Stone-Weierstrass theorem $\partial\Di\setminus F$ is dense in ${\cal M}(AP(F))$. Hence, $\widetilde p_{F_{1}F_{2}}$ and $\widetilde p_{F}$ are surjective maps.

We begin with the description of ${\cal M}(AP(F))$. 
Suppose that $F:=\{z_{1},\dots, z_{n}\}$ and $F_{i}:=F\setminus\{z_{i}\}$, $1\leq i\leq n$. Consider disjoint union 
$$
X=\bigsqcup_{1\leq i\leq n}({\cal M}(AP(\{z_{i}\}))\setminus F_{i}).
$$
Note that each component of $X$ contains $\partial\Di\setminus F$ as an open subset. By  \penalty-10000 $h_{i}:\partial\Di\setminus F\hookrightarrow {\cal M}(AP(\{z_{i}\}))\setminus F_{i}$ we denote the corresponding embeddings. Then for each $z\in\partial\Di\setminus F$ we sew together points
$h_{i}(z)$, $1\leq i\leq n$,  and identify the obtained point with $z$. As a result we obtain the quotient space $\widetilde X$ of $X$ and the ``sewing`` map $\pi:X\to\widetilde X$. We equip $\widetilde X$ with the quotient topology:
$$
U\subset\widetilde X\ \ \ {\rm is\ open}\ \ \  \Longleftrightarrow
\ \ \ \pi^{-1}(U)\subset X\ \ \ {\rm is\ open.}
$$
\begin{Lm}\label{quotient}
$\widetilde X$ is homeomorphic to ${\cal M}(AP(F))$.
\end{Lm}
{\bf Proof.} By definition each $V_{i}:=\pi({\cal M}(AP(\{z_{i}\}))\setminus F_{i})$ is an open subset of $\widetilde X$ homeomorphic to ${\cal M}(AP(\{z_{i}\}))\setminus F_{i}$. Since the latter spaces are Hausdorff, $\widetilde X$ is Hausdorff, as well. Let us cover $\partial\Di$ by closed arcs $\gamma_{1},\dots, \gamma_{n}$ such that $\gamma_{i}\cap F=\{z_{i}\}$, $1\leq i\leq n$. By $\widetilde \gamma_{i}$ we denote the closure of $\gamma_{i}$ in ${\cal M}(AP(\{z_{i}\}))$. Then 
$\widetilde \gamma_{i}$ is a compact subset of ${\cal M}(AP(\{z_{i}\}))\setminus F_{i}$ and $U_{i}:=\pi(\widetilde \gamma_{i})$ is a compact subset of $V_{i}$. It is easy to see that $\widetilde X=\cup_{1\leq i\leq n}U_{i}$. Thus $\widetilde X$ is a compact space.
Further, according to (\ref{e21}) each function from $AP(F)$ is extended continuously to $\widetilde X$ and the algebra of the extended functions separates points on $\widetilde X$. Hence by the Stone-Weierstrass theorem, $\widetilde X$ is homeomorphic to ${\cal M}(AP(F))$.\ \ \ \ \ $\Box$

As a corollary of the above construction we immediately obtain the following:

Let $F_{1}\subset F_{2}$ be finite subsets of $S$. 
Consider the commutative diagram
\begin{equation}\label{diag1}
\begin{array}{ccccc}
{\cal M}(AP(S))&\stackrel{r_{S}}{\longrightarrow}&{\cal M}(R_{S})&\stackrel{c_{S}}{\longrightarrow}&\partial\Di
\\
^{_{\widetilde p_{F_{2}}}}\!\!\downarrow&&^{_{p_{F_{2}}}}\!\!\downarrow&&||
\\
{\cal M}(AP(F_{2}))&\stackrel{r_{F_{2}}}{\longrightarrow}&{\cal M}(R_{F_{2}})&\stackrel{c_{F_{2}}}{\longrightarrow}&\partial\Di\\
^{_{\widetilde p_{F_{1}F_{2}}}}\!\!\downarrow&&^{_{p_{F_{1}F_{2}}}}\!\!\downarrow&&||
\\
{\cal M}(AP(F_{1}))&\stackrel{r_{F_{1}}}{\longrightarrow}&{\cal M}(R_{F_{1}})&\stackrel{c_{F_{1}}}{\longrightarrow}&\partial\Di.
\end{array}
\end{equation}
Here $\widetilde p_{F_{1}}:=\widetilde p_{F_{1}F_{2}}\circ \widetilde p_{F_{2}}$ and $p_{F_{1}}:=p_{F_{1}F_{2}}\circ p_{F_{2}}$. We set
$\widetilde F_{1}:=(c_{F_{1}}\circ r_{F_{1}})^{-1}(F_{1})$ and
$\widetilde S_{i}:=(c_{F_{i}}\circ r_{F_{i}})^{-1}(S)$, $i=1,2$.
Then
\begin{itemize}
\item[(A)]
$$
\widetilde p_{F_{1}F_{2}}:\widetilde p_{F_{1}F_{2}}^{-1}(\widetilde F_{1})\to
\widetilde F_{1}
$$
is a homeomorphism;
\item[(B)]
$$
{\cal M}(AP(F_{2}))\setminus\widetilde S_{2}\stackrel{\widetilde p_{F_{1}F_{2}}}{\longrightarrow} {\cal M}(AP(F_{1}))\setminus\widetilde S_{1}\stackrel{c_{F_{1}}\circ r_{F_{1}}}{\longrightarrow}\partial\Di\setminus S 
$$
are the identity maps.
\end{itemize}
From here by the definition of the inverse limit we obtain
\begin{itemize}
\item[(A1)]
$$
\widetilde p_{F_{1}}:\widetilde p_{F_{1}}^{-1}(\widetilde F_{1})\to
\widetilde F_{1}
$$
is a homeomorphism;
\item[(B1)]
$$
{\cal M}(AP(S))\setminus (c_{S}\circ r_{S})^{-1}(S) \stackrel{\widetilde p_{F_{1}}}{\longrightarrow} {\cal M}(AP(F_{1}))\setminus\widetilde S_{1}\stackrel{c_{F_{1}}\circ r_{F_{1}}}{\longrightarrow}\partial\Di\setminus S 
$$
are the identity maps.
\end{itemize}

Assume now that $F_{1}=\{z\}\subset S$. Then according to Theorem \ref{te1} (1) the set
$c_{F_{1}}^{-1}(F_{1})=c_{S}^{-1}(F_{1})$ consists of two points $\{z_{+}\}$ and $\{z_{-}\}$ identified with counterclockwise and clockwise orientations at $z$. Thus to prove (4) we must show (according to (\ref{diag1}) and
statements (A), (A1)) that each set $r_{\{z\}}^{-1}(z_{\pm})$ is homeomorphic to $b\Re$.

To this end let us recall that in the proof of Theorem \ref{te1} we had established that
${\cal M}(AP(\{z\}))$ is the inverse limit of the maximal ideal spaces ${\cal M}_{\{z\}}(s)$ of algebras $AP_{\{z\}}(s)$ of continuous functions on $\partial\Di\setminus\{z\}$ almost periodic on the open arcs
$\gamma_{t^{k}}(s)$ where $z:=e^{it}$ and $s\in (0,\pi)$.
Also, in that proof the structure of each ${\cal M}_{\{z\}}(s)$ was described. 

For every pair $0<s''<s'<\pi$ let $p_{s''s'}:{\cal M}_{\{z\}}(s'')\to {\cal M}_{\{z\}}(s')$ be the continuous surjective map dual to the embedding $i_{s''s'}: AP_{\{z\}}(s')\hookrightarrow AP_{\{z\}}(s'')$. From the proof of Theorem \ref{te1} we know that every ${\cal M}_{\{z\}}(s)$ is obtained by gluing $\partial\Di\setminus\{z\}$ with two copies of $b\Re$ where one copy (denoted $b\Re_{1}$) is obtained by gluing with $\gamma_{t^{1}}(s)$ and another one (denoted by $b\Re_{-1}$) is obtained by gluing with $\gamma_{t^{-1}}(s)$. Suppose that $\xi\in b\Re_{1}\subset {\cal M}_{\{z\}}(s'')$. Let us compute $p_{s''s'}(\xi)\in b\Re_{1}\subset
{\cal M}_{\{z\}}(s')$. Let $\{z_{\alpha}\}\subset \gamma_{t^{1}}(s'')$ be a net converging to $\xi$. This means that
the net $\{\phi_{s''}(z_{\alpha})\}\subset\Re_{-}$ converges to $\xi$ in the topology of the Bohr compactification on $b\Re$; here $\phi_{s''}$ is the map inverse to the map $\psi_{s''}:\Re_{-}\to \gamma_{t^{1}}(s'')$, $x\mapsto e^{i(t+s''e^{x})}$. Next, by the definition the net $\{\phi_{s'}(z_{\alpha})\}$ converges to 
$p_{s''s'}(\xi)$.
A straightforward computation shows that $$
\phi_{s'}(z_{\alpha})=\phi_{s''}(z_{\alpha})+\ln\left(\frac{s'}{s''}\right)\ \ \ {\rm for\ all}\ \ \ z_{\alpha}.
$$
Thus we have
\begin{equation}\label{form1}
p_{s''s'}(\xi)=\xi+\ln\left(\frac{s'}{s''}\right),\ \ \ \xi\in b\Re_{1}.
\end{equation}
Here the sum denotes the group operation on $b\Re$.
Similarly,
\begin{equation}\label{form2}
p_{s''s'}(\xi)=\xi+\ln\left(\frac{s'}{s''}\right),\ \ \ \xi\in b\Re_{-1}.
\end{equation}

Using these formulas we now prove that each $r_{\{z\}}^{-1}(z_{\pm})$ is homeomorphic to $b\Re$.
We will prove the statement for $z_{+}$ (for $z_{-}$ the argument is similar). 

For a fixed $s_{0}\in (0,\pi)$ consider the limit map
$p_{s_{0}}:{\cal M}(AP(\{z\}))\to {\cal M}_{\{z\}}(s_{0})$.
Then $p_{s_{0}}$ maps $r_{\{z\}}^{-1}(z_{+})$ into $X_{s_{0}}:=b\Re_{1}\subset {\cal M}_{\{z\}}(s_{0})$. Moreover, by definition
$r_{\{z\}}^{-1}(z_{+})$ is the inverse limit of the system $\{(X_{s''},X_{s'},p_{s''s'})\}$ where we write $X_{s}$ for $b\Re_{1}\subset {\cal M}_{\{z\}}(s)$. Since according to (\ref{form1}) every $p_{s''s'}:X_{s''}\to X_{s'}$ is a homeomorphism (even an automorphism of $b\Re$), by the definition of the inverse limit $p_{s_{0}}:r_{\{z\}}^{-1}(z_{+})\to X_{s_{0}}$ is a homeomorphism.

This completes the proof of statement (4).

(5) The required statement follows from (B1) and Theorem \ref{te2} (2).

The proof of Theorem \ref{te2} is complete.\ \ \ \ \ $\Box$
\begin{R}\label{adstat}
{\rm It is well known that the covering dimension of $b\Re$ is $\infty$ because this group is the inverse limit of compact abelian Lie groups whose dimensions tend to $\infty$.
Since $b\Re\subset {\cal M}(AP(S))$, the covering dimension of ${\cal M}(AP(S))$ is $\infty$, as well (statement (6)).
Also, statement (7) follows from the fact that $c_{S}^{-1}(S)$ is totally disconnected. Indeed, for a continuous map $\phi:T\to (c_{S}\circ r_{S})^{-1}(S)$ of a connected topological space $T$ the image of $r_{S}\circ \phi$ is a single point. This implies the required.}
\end{R}
\sect{\hspace*{-1em}. Proofs of Theorem \ref{te3} and Example \ref{notequal}} 
{\bf 4.1.} In this part we formulate and prove some auxiliary results used in the proof of the theorem.
\\
{\bf Notation.} {\em
Let $z_{0}\in \partial\Di$ and $U_{z_{0}}$ be the intersection of an open disk of radius $\leq 1$ centered at $z_{0}$ with $\overline{\Di}\setminus\{z_{0}\}$. We call such $U_{z_0}$ a \textit{circular neighbourhood} of $z_0$.}

Next, we define almost periodic continuous functions on a circular neighbourhood $U_{z_{0}}$ of $z_{0}$ holomorphic in its interior points as follows.

Let $\phi_{z_{0}}:\Di\to\H_{+}$,
\begin{equation}\label{fi}
\phi_{z_{0}}(z):=\frac{2i(z_{0}-z)}{z_{0}+z},\ \ \ z\in\Di,
\end{equation}
be a conformal map of $\Di$ onto the upper half-plane $\H_{+}$. Then
$\phi_{z_{0}}$ is also continuous on $\partial\Di\setminus\{-z_{0}\}$ and maps it diffeomorphically onto $\Re$ (the boundary of $\H_{+}$) so that $\phi_{z_{0}}(z_{0})=0$.
Let $\Sigma_{0}$ be the interior of the strip $\Sigma:=\{z\in\Co\ :\ {\rm Im}\ z\in [0,\pi]\}$. Consider the conformal map ${\rm Log}:\H_{+}\to\Sigma_{0}$,
$z\mapsto {\rm Log}(z):=\ln|z|+i{\rm Arg}(z)$, where ${\rm Arg}:\Co\setminus\Re_{-}\to (-\pi,\pi)$ is the principal branch of the multi-function ${\rm arg}$, the argument of a complex number. The function ${\rm Log}$ is extended to a homeomorphism of $\overline{\H}_{+}\setminus\{0\}$ onto $\Sigma$; here $\overline{\H}_{+}$ stands for the closure of $\H_{+}$.

By $AP_{\cal C}(\Sigma)$ we denote the algebra of uniformly continuous almost periodic functions on $\Sigma$ (i.e., they are almost periodic on any line parallel to the $x$-axis). Clearly we have $AP_{\cal O}(\Sigma)\subset AP_{\cal C}(\Sigma)$. Then according to Theorem \ref{coron1} (the corona theorem for $AP_{\cal O}(\Sigma)$), the maximal ideal space ${\cal M}(AP_{\cal C}(\Sigma))$ of the algebra $AP_{\cal C}(\Sigma)$ is homeomorphic to ${\cal M}(AP_{\cal O}(\Sigma))$. In what follows we identify these spaces.
\begin{D}\label{apdisc}
We say that $f:U_{z_{0}}\to\Co$ is a (continuous) almost periodic function if there is a function $\widehat f\in AP_{{\cal C}}(\Sigma)$ such that 
$$
f(z):=\widehat f({\rm Log}(\phi_{z_{0}}(z)))\ \ \ {\rm for\ all}\ \ \ z\in U_{z_{0}}.
$$
If such $\widehat f\in AP_{{\cal O}}(\Sigma)$, then $f$ is called a holomorphic almost periodic function.
\end{D}

Suppose that $z_{0}=e^{it_{0}}$. For $s\in (0,\pi)$ we set $\gamma_{1}(z_0,s):={\rm Log}(\phi_{z_{0}}(\gamma_{t_{0}^{1}}(s)))\subset\Re$ and 
$\gamma_{-1}(z_0,s):={\rm Log}(\phi_{z_{0}}(\gamma_{t_{0}^{-1}}(s)))\subset \Re+i\pi$. 
\begin{Lm}\label{le1}
Let $f\in AP(\{-z_{0},z_{0}\})$. We set $f_{k}:=f|_{\gamma_{t_{0}^{k}}(\pi)}$ and
consider the functions $h_{k}=f_{k} \circ \varphi_{z_0}^{-1} \circ \rm Log^{-1}$ on $\gamma_{k}(z_0,s)$, $k\in\{-1,1\}$. Then for any $\epsilon>0$ there are $s_{\epsilon}\in (0,s)$ and almost periodic functions $h_{1}'$ on $\Re$, and $h_{-1}'$ on $\Re+i\pi$ such that 
\begin{equation}
\sup_{z\in \gamma_{k}(z_{0},s_{\epsilon})}|h_{k}(z)-h_{k}'(z)|<\epsilon\ \ \ {\rm for\ each}\ \ \ k\in\{-1,1\}.
\end{equation}
\end{Lm}
{\bf Proof.} We prove the result for $f_{1}$ only. The proof for $f_{-1}$ is similar. According to Theorem \ref{te1} it suffices to prove the lemma for $f_{1}=g_{t_{0}}$ or $f_{1}=e^{i\lambda\log_{t_{0}^{1}}}$, $\lambda\in\Re$. In the first case we can choose a sufficiently small $s_{\epsilon}$ such that on $\gamma_{t_{0}^{1}}(s_{\epsilon})$ the function $g_{t_{0}}$ is uniformly approximated with an error $<\epsilon$ by a constant function. Then as $h_{1}'$ we can choose the corresponding constant function on $\gamma_{1}(z_{0},s_{\epsilon})$. In the second case, by definition,
$$
h_{1}(x)=e^{i\lambda\ln\left( {\rm Arg}\left(\frac{2i-e^{x}}{2i+e^{x}}\right)\right)}=
e^{i\lambda\ln\left(\frac{4e^{x}}{4-e^{2x}}+o(e^{3x})\right)} =e^{i\lambda(x+o(e^{x}))}\ \ \ {\rm as}\ \ \ x\to-\infty.
$$
From here for a sufficiently small $s_{\epsilon}$ we have
$$
|h_{1}(x)-e^{i\lambda x}|<\epsilon\ \ \ {\rm for\ all}\ \ \ x\in \gamma_{1}(z_{0},s_{\epsilon}).
\ \ \ \ \ \Box
$$

We also use the following well-known result.
\begin{Lm}\label{pois}
Suppose that $f_{1}$ and $f_{2}$ are continuous almost periodic functions on $\Re$ and $\Re+i\pi$, respectively. Then there exists a function $F \in AP_{\cal C}(\Sigma)$ harmonic in $\Sigma_{0}$ whose boundary values are $f_{1}$ and $f_{2}$.
\end{Lm}
{\bf Proof.} Let $F$ be a function harmonic in $\Sigma_{0}$ with boundary values $f_{1}$ and $f_{2}$. Since $f_{1}$ and $f_{2}$ are almost periodic, for any $\epsilon>0$ there exists $l(\epsilon)>0$ such that every interval $[t_{0},t_{0}+l(\epsilon)]$ contains a common $\epsilon$-period of $f_{1}$ and $f_{2}$, say, $\tau_{\epsilon}$, see, e.g., [LZ]. Thus
\begin{displaymath}
\sup_{x \in \Re}|f_{1}(x+\tau_{\epsilon})-f_{1}(x)|<\epsilon\ \ \ {\rm and }
\ \ \ \sup_{x \in \Re}|f_{2}(x+i\pi+\tau_{\epsilon})-f_2(x+i\pi)|<\epsilon.
\end{displaymath}
Now, by the maximum principle for harmonic functions
\begin{displaymath}
\sup_{x \in \Re}|F(x+iy+\tau_\epsilon)-F(x+iy)|<\epsilon\ \ \ {\rm for\ each}\ \ \ y\in [0,\pi],
\end{displaymath}
that is, $F$ is almost periodic on every line $\Re+iy$, $y \in [0,\pi]$.
\ \ \ \ \ $\Box$

Let ${\cal A}(U_{z_{0}})$ be the algebra of continuous functions on $\Di$ almost periodic on the circular neighbourhood $U_{z_{0}}$ of $z_{0}$. By 
${\cal A}_{z_{0}}$ we denote the uniform closure of the algebra generated by all ${\cal A}(U_{z_{0}})$ and
by ${\cal M}_{z_{0}}$ the closure of $\Di$ in the maximal ideal space of ${\cal A}_{z_{0}}$. Since the algebra ${\cal A}_{z_{0}}$ is self-adjoint, by the Stone-Weierstrass theorem ${\cal M}_{z_0}$ coincides with the maximal ideal space of ${\cal A}_{z_{0}}$. Next, let $p_{z_{0}}:{\cal M}_{z_{0}}\to\overline\Di$ be the continuous surjective map dual to the natural embedding $C(\overline\Di)\hookrightarrow {\cal A}_{z_{0}}$.
Then we have
\begin{Lm}\label{le3}
\begin{itemize}
\item[(a)]
For every neighbourhood $U$ of the compact set $F_{z_{0}}:=p_{z_{0}}^{-1}(z_{0})$ there is a circular neighbourhood $U_{z_{0}}$ of $z_{0}$ such that  $U_{z_{0}}\cap\Di\subset U\cap \Di$.
\item[(b)]
$F_{z_{0}}$ is homeomorphic to ${\cal M}({\cal}{AP}_{\cal O}(\Sigma))$.
\item[(c)]
Each function $f\in AP(S)\cap H^{\infty}$ belongs to the algebra $\displaystyle \bigcap_{z\in\partial\Di}{\cal A}_{z}$.
\end{itemize}
\end{Lm}
{\bf Proof.} (a), (b). Since the algebra ${\cal A}(U_{z_{0}})$ is self-adjoint, $\Di$ is dense in its maximal ideal space ${\cal M}(U_{z_{0}})$. Then ${\cal M}_{z_{0}}$ is the inverse limit of the compact spaces ${\cal M}(U_{z_{0}})$ (because ${\cal A}_{z_{0}}$ is the uniform closure of the algebra generated by algebras ${\cal A}(U_{z_{0}})$), see, e.g, [R].  For $U_{z_{0}}\subset V_{z_{0}}$ by
$p_{U_{z_{0}}V_{z_{0}}}:{\cal M}(U_{z_{0}})\to {\cal M}(V_{z_{0}})$ we denote the maps in this limit system and by $p_{U_{z_{0}}}:{\cal M}_{z_{0}}\to {\cal M}(U_{z_{0}})$ the corresponding (continuous and surjective) limit maps.
Then, by the definition of the inverse limit, the base of topology on ${\cal M}_{z_{0}}$ consists of the sets
$p_{U_{z_{0}}}^{-1}(U)$ where $U\subset {\cal M}(U_{z_{0}})$ is open and  $U_{z_{0}}$ is a circular neighbourhood of $z_{0}$.
In particular, since $F_{z_{0}}$ is a compact set, for a neighbourhood $U$ of $F_{z_{0}}$ there is a circular neighbourhood $\widetilde U_{z_{0}}$ of $z_{0}$ and a neighbourhood $\widetilde U\subset {\cal M}(\widetilde U_{z_{0}})$ of $F(\widetilde U_{z_{0}}):=p_{\widetilde U_{z_{0}}}(F_{z_{0}})$ such that $p_{\widetilde U_{z_{0}}}^{-1}(\widetilde U)\subset U$. Recall that $\Di$ is a dense subset of ${\cal M}(\widetilde U_{z_{0}})$ and ${\cal M}_{z_{0}}$. Also,
$p_{\widetilde U_{z_{0}}}^{-1}(\Di)$ contains $\Di\subset {\cal M}_{z_{0}}$. Thus in order to prove (a) it suffices to show that there is $U_{z_{0}}\subset \widetilde U_{z_{0}}$ such that $U_{z_{0}}\cap\Di\subset\widetilde U\cap\Di$.

First, let us study the structure of ${\cal M}(\widetilde U_{z_{0}})$. 
Let ${\cal A}^{*}(\widetilde U_{z_{0}})$ be the pullback to $\Sigma$ by means of the map $({\rm Log}\circ\phi_{z_{0}})^{-1}$ of the algebra ${\cal A}(\widetilde U_{z_{0}})$.
Then ${\cal A}^{*}(\widetilde U_{z_{0}})$ consists of continuous functions on $\Sigma_{0}$ such that on $({\rm Log}\circ\phi_{z_{0}})(\widetilde U_{z_{0}})$ they are restrictions of almost periodic functions on $\Sigma$. Since ${\cal A}^{*}(\widetilde U_{z_{0}})$ is isomorphic to ${\cal A}(\widetilde U_{z_{0}})$ we can naturally identify the maximal ideal spaces of these algebras.

Further, observe that there is $T<0$ such that $({\rm Log}\circ\phi_{z_{0}})(\widetilde U_{z_{0}})$ contains the subset $\Sigma_{T}:=\{z\in\Sigma\ :\ {\rm Re}\ z\leq T\}$ of the strip $\Sigma$. By the definition of the topology on ${\cal M}(AP_{\cal O}(\Sigma))$ (see section 2) $\Sigma_{T}$ is dense in ${\cal M}(AP_{\cal O}(\Sigma))$. Hence,
the space ${\cal M}(\widetilde U_{z_{0}})$ contains ${\cal M}(AP_{{\cal C}}(\Sigma)) (={\cal M}(AP_{{\cal O}}(\Sigma)))$. Let $K$ be the intersection of the closures of $\widetilde U_{z_{0}}$ and $\Di\setminus \widetilde U_{z_{0}}$ in $\Co$. Then
$K':={\rm Log}\circ\phi_{z_{0}}(K)$ is a compact subset of $\Sigma$. In particular, $AP_{\cal C}(\Sigma)|_{K'}=C(K')$. This implies (by the Tietze extension theorem) that every bounded continuous function on $(\Di\setminus \widetilde U_{z_{0}})\cup K$ can be extended to a function from ${\cal A}(\widetilde U_{z_{0}})$ with the same supremum norm.
Now, we can describe explicitly ${\cal M}(\widetilde U_{z_{0}})$ as follows (cf. the proof of Theorem \ref{te1} for a similar construction).

Let $M$ be the maximal ideal space of the algebra of bounded continuous functions on $(\Di\setminus \widetilde U_{z_{0}})\cup K$. We identify $K\subset M$  with $K'\subset {\cal M}(AP_{\cal O}(\Sigma))$ by means of ${\rm Log}\circ\phi_{z_{0}}$. Then the quotient space of $M\sqcup {\cal M}(AP_{\cal O}(\Sigma))$ under this identification is homeomorphic to ${\cal M}(\widetilde U_{z_{0}})$.

Next, by definition the fibre $F_{z_{0}}$ over $z_{0}$ consists of limit points in ${\cal M}_{z_{0}}$ of all nets converging to $\{z_{0}\}$ inside $\Di$. This and the above construction of ${\cal M}(\widetilde U_{z_{0}})$ show that
$F(\widetilde U_{z_{0}})$ is homeomorphic to ${\cal M}(AP_{{\cal O}}(\Sigma))$. Moreover, $F_{z_{0}}$ is the inverse limit of compact sets $F(\widetilde U_{z_{0}})$ where the limit system is determined by the maps $p_{U_{z_{0}}V_{z_{0}}}|_{F(U_{z_{0}})}$ for $U_{z_{0}}\subset V_{z_{0}}$. 
Let us show that the maps $p_{U_{z_{0}}V_{z_{0}}}|_{F(U_{z_{0}})}:F(U_{z_{0}})\to F(V_{z_{0}})$ are homeomorphisms. Indeed, let $\{z_{\alpha}\}\subset\ U_{z_{0}}$ be a net converging to a point $\xi\in F(U_{z_{0}})$. Since $U_{z_{0}}\hookrightarrow V_{z_{0}}$, in our definitions of ${\cal M}(U_{z_{0}})$ and ${\cal M}(V_{z_{0}})$ the net $\{z_{\alpha}\}$ converges to the same point $\xi\in F(V_{z_{0}})$ which gives the required statement. Since all maps $p_{U_{z_{0}}V_{z_{0}}}|_{F(U_{z_{0}})}$ are homeomorphisms, 
by the definition of the inverse limit the map $p_{\widetilde U_{z_{0}}}|_{F_{z_{0}}}:F_{z_{0}}\to F(\widetilde U_{z_{0}})$,
is also a homeomorphism. This completes the proof of (b).

Now, observe that in our model of ${\cal M}(\widetilde U_{z_{0}})$ the intersection of $\widetilde U$ with $\Di$ contains $\widetilde U_{z_{0}}\cap\Di$, this completes the proof of (a).

(c) Fix a point $z_{*}\in\partial\Di$. We must show that every $f\in AP(S)\cap H^{\infty}$ belongs to ${\cal A}_{z_{*}}$. According to (\ref{e21}) and Lemma \ref{le1} each $f\in AP(\partial\Di)\cap H^{\infty}$
can be approximated locally on open arcs of the form $\gamma_{t^{k}}(s)$, $k\in\{-1,1\}$, $s\in (0,\pi)$, $z:=e^{it}\in\partial\Di$, by 
pullbacks of almost periodic functions on the boundary of $\Sigma$. Using compactness of $\partial\Di$, for a given $\epsilon>0$  we can find a finitely many points $z_{1},\dots, z_{n}\in \partial\Di$, arcs $\gamma_{t_{l}^{k}}(s_{l})$, $k\in\{-1,1\}$, $s_{l}\in (0,\pi)$, $z_{l}:=e^{it_{l}}$, and functions $f^{l}:\partial\Di\setminus\{-z_{l},z_{l}\}\to\Co$ which are pullbacks of almost periodic functions on $\partial\Sigma$ by means of  ${\rm Log}\circ\phi_{z_{l}}$, $1\leq l\leq n$, such that  
$$
\begin{array}{c}
\displaystyle
\partial\Di\setminus\{z_{1},\dots, z_{n}\}=\bigcup_{1\leq l\leq n}(\gamma_{t_{l}^{-1}}\cup \gamma_{t_{l}^{1}})\ \ \ {\rm and\ for\ all}\ \ \ 1\leq l\leq n\\
\\
\displaystyle
{\rm ess}\!\!\!\!\!\!\!\!\!\!\!\!\!\!\sup_{z \in \gamma_{t_{l}^{1}}(s_{l})\cup \gamma_{t_{l}^{-1}}(s_{l})}|f(z)-f^l(z)|<\epsilon.
\end{array}
$$
Without loss of generality we may assume that  $z_{*}\in\{z_{1},\dots, z_{n}\}$.
Next, set  $V_{l}:=\gamma_{t_{l}^{-1}}\cup \gamma_{t_{l}^{1}}\cup\{z_{l}\}$. Then 
$(V_{l})_{l=1}^n$ is a finite open cover of $\partial\Di$. Let $\{\rho_{l}\}_{l=1}^n$ be a smooth partition of unity subordinate to this cover such that $\rho_{l}(z_{l})=1$.
Consider the functions $f_{l}$ on $\partial\Di\setminus\{z_{l}\}$ defined by the 
formulas
$$
f_{l}:=\rho_{l}f^{l},\ \ \ 1\leq l\leq n.
$$
Let $F$ be the harmonic function on $\Di$ such that
$F|_{\partial\Di}=\sum_{1\leq l\leq n}f_{l}$.
Then 
$$
||f-F||_{L^{\infty}(\Di)}<\epsilon.
$$
We prove that $F\in {\cal A}_{z_{*}}$. Since $\epsilon>0$ is arbitrary, this will complete the proof of (c).

By $F_{l,1}$ and $F_{l,2}$ we denote the harmonic functions on $\Di$ with the boundary values $f_{l}$ and $f^{l}-f_{l}$, respectively. Thus $F_{l}:=F_{l,1}+F_{l,2}$ is the harmonic function with the boundary values $f^{l}$. According to Lemma \ref{pois}, every $F_{l}$ is almost periodic on $\overline{\Di}\setminus\{\pm z_{l}\}$. Thus if $z_{l}=z_{*}$, then $F_{l}\in {\cal A}_{z_{*}}$. If $z_{l}\neq z_{*}$, then $F_{l}$ is continuous at $z_{*}$ and so by the definition of ${\cal A}_{z_{*}}$ the function $F_{l}\in {\cal A}_{z_{*}}$, as well.
Further, for a point $z_{l}$ distinct from $z_{*}$ the function $F_{l,1}$  can be extended continuously in an open disk centered at $z_{*}$ (because the support of $f_{l}$ does not contain $z_{*}$). Hence,
such $F_{l,1}\in {\cal A}_{z_{*}}$. Assume now that $z_{l}=z_{*}$ for some $l$.
Then the function $F_{l,2}$ can be extended continuously in an open disk centered at $z_{*}$ (because the support of $f^{l}-f_{l}$ does not contain $z_{*}$). Thus $F_{l,2}\in {\cal A}_{z_{*}}$ and in this case
$F_{l,1}:=F_{l}-F_{l,2}\in {\cal A}_{z_{*}}$, as well. Since
$F:=\sum_{1\leq l\leq n}F_{l,1}$, combining the above considered cases we obtain that $F\in {\cal A}_{z_{*}}$. 

This completes the proof of Lemma \ref{le3}.\ \ \ \ \ $\Box$

\begin{Th}\label{appr1}
Let $f\in AP(S)\cap H^{\infty}$. Then for each $z_{0}\in \partial\Di$ and any $\epsilon>0$ there is a circular neighbourhood $U_{z_{0}}:=U_{z_{0}}(f,\epsilon)$ of $z_{0}$ and a holomorphic almost periodic function $f_{z_{0}}$ on $U_{z_{0}}$ such that
$$
\sup_{z\in U_{z_{0}}\cap\Di}|f(z)-f_{z_{0}}(z)|<\epsilon.
$$
\end{Th}
{\bf Proof.} Fix a point $z_{0}\in\partial\Di$.
According to Lemma \ref{le3}(c) the function \penalty-10000 $f\in AP(S)\cap H^{\infty}$ belongs to ${\cal A}_{z_{0}}$ and so it is extended by means of the Gelfand transform to a continuous function
$\widehat f$ on ${\cal M}_{z_{0}}$. We use the description of ${\cal M}_{z_{0}}$ presented in the proof of Lemma \ref{le3}. Recall that in that construction the fibre $F_{z_{0}}\subset {\cal M}_{z_{0}}$ over $z_{0}$  is naturally identified with ${\cal M}(AP_{\cal O}(\Sigma))$.
Then we have (see Definition \ref{holom1}).
\begin{Lm}\label{holext}
The function $\widehat f|_{F_{z_{0}}}$ is holomorphic.
\end{Lm}
{\bf Proof.} In the proof we use the results of section 2. 
Let us consider the map $i_{\xi}:\Sigma_{0}\to {\cal M}(AP_{\cal O})(\Sigma)$, $\xi\in b\Z$. We must show that $\widehat f\circ i_{\xi}$ is holomorphic.

To this end we transfer the function $f$ by means of the map $({\rm Log}\circ\phi_{z_{0}})^{-1}$ to $\Sigma_{0}$ and by $\widetilde f$ denote the pulled back function. Fix a point $\eta\in i_{\xi}(\Sigma_{0})$, say, $\eta:=i_{\xi}(w)$, $w\in\Sigma_{0}$.
Then there is a straight line $\Re+iy\subset\Sigma_{0}$, $y\in (0,\pi)$, and
a net $\{z_{\alpha}\}\subset\Re+iy$  which forms an infinite discrete subset of $\Re+iy$ such that $\{z_{\alpha}\}$ converges to $\eta$ in the topology of ${\cal M}(AP_{\cal O})(\Sigma)$. 
Let $B\subset\Co$ be an open disk such that $\{\Re+iy\}+ B\subset\Sigma_{0}$. By the definition of $i_{\xi}$, for each $z\in B$ the net $\{z_{\alpha}+z\}$ converges in
${\cal M}(AP_{\cal O})(\Sigma)$ to $i_{\xi}(w+z)$. Also, by the definition of ${\cal M}_{z_{0}}$ we have
$$
\lim_{\alpha}\widetilde f(z_{\alpha}+z)=(\widehat f\circ i_{\xi})(w+z),\ \ \ z\in B.
$$
But the holomorphic functions $\widetilde f_{\alpha}(z):=\widetilde f(z_{\alpha}+z)$ form a normal family on $B$. Therefore using an argument similar to that of the proof of Lemma  \ref{holom3} (1) we obtain that $\widehat f\circ i_{\xi}|_{B}$ is holomorphic.
Since $\xi$ and $\eta$ are arbitrary, the latter implies that $\widehat f|_{F_{z_{0}}}$ is holomorphic.\ \ \ \ \ $\Box$

Now by Lemma \ref{holom3} we obtain that there is a function $\widetilde f_{z_{0}}\in AP_{\cal O}(\Sigma)$ whose extension to ${\cal M}(AP_{\cal O}(\Sigma))$ coincides with $\widehat f|_{F_{z_{0}}}$. Let us consider the function $f_{z_{0}}\in {\cal A}_{z_{0}}$ whose pullback to $\Sigma$ by means of $({\rm Log}\circ\phi_{z_{0}})^{-1}$ coincides with $\widetilde f_{z_{0}}$. Then by the definition of the topology of ${\cal M}_{z_{0}}$ the extension $\widehat f_{z_{0}}$ of $f_{z_{0}}$ to ${\cal M}_{z_{0}}$ satisfies $\widehat f_{z_{0}}|_{F_{z_{0}}}=\widehat f|_{F_{z_{0}}}$. Since $F_{z_{0}}$ is a compact set, the latter implies that there is a neighbourhood $U$ of $F_{z_{0}}$ in ${\cal M}_{z_{0}}$ such that $|\widehat f_{z_{0}}(x)-\widehat f(x)|<\epsilon$ for all $x\in U$.
Finally, by Lemma \ref{le3}(a) there is a circular neighbourhood $U_{z_{0}}$ such that
$U_{z_{0}}\cap\Di\subset U\cap\Di$. Thus $|f_{z_{0}}(z)-f(z)|<\epsilon$ for
all $z\in U_{z_{0}}\cap\Di$.\ \ \ \ \ $\Box$\\
\\
{\bf 4.2. Proof of Theorem \ref{te3}.} We must show that $A_{S}=AP(S)\cap H^{\infty}$. We split the proof into several parts.
First we prove the following statement.
\begin{Lm}\label{finite}
$AP(S)\cap H^{\infty}$ is the uniform closure of the algebra generated by all possible subalgebras $AP(F)\cap H^{\infty}$ with finite $F\subset S$. 
\end{Lm}
Then we will prove that $A_{F}=AP(F)\cap H^{\infty}$ for every finite subset $F\subset \partial\Di$. Together with the above lemma and the fact that $A_{S}$ is the uniform closure of the algebra generated by all possible subalgebras $A_{F}$ with finite $F\subset S$, this will clearly complete the proof of the theorem.\\
{\bf Proof.} According to Theorem \ref{appr1}, for a given $f\in AP(S)\cap H^{\infty}$ we can find finitely many points $z_{1},\dots, z_{n}$, circular neighbourhoods $U_{z_{1}},\dots, U_{z_{n}}$ and holomorphic almost periodic functions $f_{1},\dots, f_{n}$ defined on $U_{z_{1}},\dots, U_{z_{n}}$, respectively, such that $(U_{z_{i}})_{1\leq i\leq n}$ forms an open cover of $\partial\Di\setminus\{z_{1},\dots, z_{n}\}$ and
$$
\max_{i}||f|_{U_{z_{i}}}-f_{i}||_{L^{\infty}(U_{z_{i}})}<\epsilon.
$$
Since the discontinuities of $f|_{\partial\Di}$ belong to the closed set $S$, each function $f_{i}$ with $z_{i}\not\in S$ can be chosen also to be continuous on the closure $\overline{U}_{z_{i}}$.

Further, form a cocycle $\{c_{ij}\}$ on the intersections of sets from the above cover by the formula
$$
c_{ij}(z):=f_{i}(z)-f_{j}(z),\ \ \ z\in U_{z_{i}}\cap U_{z_{j}}.
$$
Diminishing, if necessary, the sets of the above cover we 
may assume without loss of generality that all $U_{z_{i}}\cap U_{z_{j}}$, $i\neq j$, do not contain points $z_{1},\dots, z_{n}$. Then each $U_{z_{i}}\cap U_{z_{j}}$, $i\neq j$, is a compact subset of $\overline\Di$ and the corresponding $c_{ij}$ are continuous and holomorphic in interior points of $U_{z_{i}}\cap U_{z_{j}}$. 

Let $\{\rho_{i}\}$ be a smooth partition of unity subordinate to the cover $(U_{z_{i}})_{1\leq i\leq n}$. We can choose every $\rho_{i}$ so that it is the restriction to $\overline{U}_{z_{i}}$ of a $C^{\infty}$-function on $\Co$ and $\rho_{i}(z_{i})=1$. As usual, we resolve the cocycle $\{c_{ij}\}$ using this partition of unity by the formulas
\begin{equation}\label{resol}
\widetilde f_{j}(z)=\sum_{k=1}^{n}\rho_{k}(z)c_{jk}(z),\ \ \ z\in \overline{U}_{z_{j}}.
\end{equation}
Hence,
$$
c_{ij}(z):=\widetilde f_{i}(z)-\widetilde f_{j}(z),\ \ \ z\in U_{z_{i}}\cap U_{z_{j}}.
$$
In particular, since $c_{ij}$ are holomorphic in $\Di\cap U_{z_{i}}\cap U_{z_{j}}$, the formula
$$
h(z):=\frac{\partial\widetilde f_{i}(z)}{\partial\overline z},\ \ \ z\in U_{z_{i}}\cap\Di,
$$
determines a smooth bounded function in an open annulus 
$A\subset\cup_{i=1}^{n}\overline U_{z_{i}}$ with the outer boundary $\partial\Di$. Also, by our choice of the partition of unity,  $h$ is extended continuously to the closure $\overline{A}$ of $A$. 

Let us consider a function $H$ determined by the formula
\begin{equation}\label{solution}
H(z)=\frac{1}{2\pi i}\int\int_{\zeta\in A}\frac{h(\zeta)}{\zeta-z}\ \!d\zeta\wedge d\overline\zeta,\ \ \ z\in\overline{A}.
\end{equation}
Passing in (\ref{solution}) to polar coordinates with the origin at $z$, we easily obtain
\begin{equation}\label{solution1}
\sup_{z\in A}|H(z)|\leq C w(A)\sup_{z\in A}|h(z)|
\end{equation}
where $w(A)$ is the width of $A$ and $C>0$ is an absolute constant. Moreover,
$H\in C(\overline A)$ and $\partial H/\partial\overline{z}=h$ in $A$, see, e.g., [G, Chapter VIII]. Let us replace $A$ by a similar annulus of a smaller width such that for this new $A$ 
$$
\sup_{z\in A}|H(z)|<\epsilon.
$$
Now we set
$$
c_{i}(z):=\widetilde f_{i}(z)-H(z),\ \ \ z\in\overline U_{z_{i}}\cap\overline A.
$$
Then each $c_{i}$ is continuous on $U_{z_{i}}\cap\overline A$ holomorphic in interior points of this set and
$$
c_{i}(z)-c_{j}(z)=c_{ij}(z), \ \ \ z\in U_{z_{i}}\cap U_{z_{j}}.
$$
Since every $|c_{ij}(z)|<2\epsilon$ for all $z\in U_{z_{i}}\cap U_{z_{j}}$, 
$$
|c_{i}(z)|<3\epsilon,\ \ \ z\in U_{z_{i}}\cap U_{z_{j}}.
$$

Let us determine a global function $f_{\epsilon}$ on $\overline A\setminus\{z_{1},\dots, z_{n}\}$ by the formulas 
$$
f_{\epsilon}(z):=f_{i}(z)-c_{i}(z),\ \ \ z\in U_{i}\cap\overline{A}.
$$ 
Since for $z_{i}\not\in S$, the function $f_{i}$ is continuous on $\overline{U}_{z_{i}}$, from the above construction we obtain that $f_{\epsilon}\in H^{\infty}(A)\cap AP(F)$ where  $F:=\{z_{1},\dots, z_{n}\}\cap S$.  Also,
$$
||f-f_{\epsilon}||_{L^{\infty}(A)}< 4\epsilon.
$$
Let $B$ be an open disk centered at 0 whose intersection with $A$ is an annulus of width $<\epsilon$.  Consider the cocycle $c$ on $B\cap A$ defined by
$$
c(z)=f(z)-f_{\epsilon}(z),\ \ \ z\in B\cap A.
$$
By definition, $|c(z)|\leq 4\epsilon$ for all $z\in B\cap A$. 
Let $A'$ be the open annulus with the interior boundary coinciding with the interior boundary of $A$ and with the outer boundary $\{z\in\Co\ :\ |z|=2\}$.
Then $A'\cap B=A\cap B$. Consider a smooth partition of unity subordinate to the cover $\{A',B\}$ of $\Di$ which consists of smooth radial functions $\rho_{1}$ and $\rho_{2}$ such that 
$$
\max_{i}||\nabla\rho_{i}||_{L^{\infty}(\Co)}\leq \widetilde C w(B\cap A)<\widetilde C\epsilon
$$
for some absolute constant $\widetilde C>0$.
Then using arguments similar to the above based on versions of (\ref{resol}), (\ref{solution}) and (\ref{solution1}) for the cocycle $c$ and the partition of unity $\{\rho_{1},\rho_{2}\}$, we can find holomorphic functions $\overline{c}_{1}$
on $B$ and $\overline{c}_{2}$ on $A$ continuous on the corresponding boundaries such
that
$$
\begin{array}{c}
\displaystyle
\overline{c}_{1}(z)-\overline{c}_{2}(z)=c(z),\ \ \ z\in B\cap A,\ \ \
{\rm and}\\
\\
\displaystyle
\max\{||\overline{c}_{1}||_{H^{\infty}(B)}, ||\overline{c}_{2}||_{H^{\infty}(A)}\}\leq \overline{C}||c||_{H^{\infty}(A\cap B)}
\end{array}
$$ 
where $\overline{C}>0$ is an absolute constant.
Finally, define 
$$
F_{\epsilon}(z):=\left\{
\begin{array}{ccc}
f(z)-\overline{c}_{1}(z),&{\rm if}&z\in B,\\
\\
f_{\epsilon}(z)-\overline{c}_{2}(z),&{\rm if}&z\in\Di\setminus B. 
\end{array}
\right.
$$
Clearly we have
$$
||f-F_{\epsilon}||_{L^{\infty}(\Di)}<c\epsilon.
$$
for some absolute constant $c>0$, and $F_{\epsilon}\in AP(F)\cap H^{\infty}$, where $F:=\{z_{1},\dots, z_{n}\}\cap S$. Since $\epsilon$ is arbitrary, this completes the proof of the lemma.\ \ \ \ \ $\Box$

As the next step of the proof let
us establish Theorem \ref{te3} for $AP(F)\cap H^{\infty}$ with a finite set $F\subset\partial S$, say, $F=\{z_{1},\dots, z_{n}\}$.
\begin{Lm}\label{fin1}
$$
A_{F}=AP(F)\cap H^{\infty}.
$$
\end{Lm}
{\bf Proof.} Let us show that $AP(F)\cap H^{\infty}\subset A_{F}$.
First suppose that $F$ contains at least two points.

Let $\psi_{1}:\partial\Di\to\partial\Di$ be the restriction to the boundary of a M\"{o}bius transformation of $\Di$ that maps $-z_{1}$ to a point of $F$ distinct from $z_{1}$ and preserves $z_{1}$. Then by the definition of M\"{o}bius transformations $\psi_{1}$ is a $C^{1}$ diffeomorphism of $\partial\Di$. In particular, by Corollary \ref{cor1} for a given $f\in AP(F)\cap H^{\infty}$ the function $f\circ\psi_{1}\in AP(F_{1})\cap H^{\infty}$ where $F_{1}:=\psi_{1}^{-1}(F)$. Since $z_{1}\in F_{1}$,
as in the proof of Theorem \ref{appr1} we can find an almost periodic holomorphic function $g_{1}$ on $\overline\Di\setminus\{\pm z_{1}\}$ such that $f\circ\psi_{1}-g_{1}$ is continuous and equals $0$ at $z_{1}$. We set
$$
\widetilde g_{1}(z):=\frac{g(z)(z+z_{1})}{2z_{1}},\ \ \ z\in\overline\Di\setminus\{z_{1}\}.
$$
Then $\widetilde g_{1}$ has a discontinuity at $z_{1}$ only.
Let us show that $\widetilde g_{1}\in A_{\{-z_{1},z_{1}\}}$. Indeed, 
by the definition of $g_{1}$ the function  $g_{1}\circ\phi_{z_{1}}^{-1}\circ{\rm Log}^{-1}$ belongs to $AP_{{\cal O}}(\Sigma)$. Therefore, it can be uniformly approximated on $\Sigma$ by 
polynomials in variables $e^{i\lambda z}$, $\lambda\in\Re$, see, e.g., [JT].
In turn, $g_{1}$ can be uniformly approximated on $\overline\Di\setminus\{\pm z_{1}\}$ by complex polynomials in variables $e^{i\lambda{\rm Log}\circ\phi_{1}}$.
Now for $z\in \partial\Di$ we have 
$$
{\rm Im}\{({\rm Log}\circ\phi_{z_{1}})(z)\}:=\left\{
\begin{array}{ccc}
0,&{\rm if}&0\leq {\rm Arg}(z/z_{1})<\pi,\\
\\
\pi,&{\rm if}&0\leq {\rm Arg}(z_{1}/z)<\pi.
\end{array}
\right.
$$
This implies that every function $e^{i\lambda{\rm Log}\circ\phi_{1}}$, $\lambda\in\Re$,  belongs to $A_{\{-z_{1},z_{1}\}}$ and so $g_{1}\in A_{\{-z_{1},z_{1}\}}$, as well. Since $(z+z_{1})/2z_{1}\in A_{0}$, the function $\widetilde g_{1}\in A_{\{-z_{1},z_{1}\}}$ by definition.
Thus the function $h_{1}:=\widetilde g_{1}\circ\psi_{1}^{-1}$ belongs to $A_{F}$, is continuous outside $z_{1}$ and 
$f_{1}:=f-h_{1}\in AP(F^{1})\cap H^{\infty}$, where $F^{1}:=F\setminus\{z_{1}\}$. Further, using similar arguments we can find $h_{2}\in A_{F}$, continuous outside $z_{2}$ such that $f_{2}:=f_{1}-h_{2}\in AP(F^{2})\cap H^{\infty}$ where $F^{2}:=F^{1}\setminus\{z_{2}\}$ etc. After $n$ steps we obtain
functions $h_{1},\dots, h_{n}\in A_{F}$ such that $h_{k}$ is continuous outside $z_{k}$, $1\leq k\leq n$, and the function  
\begin{equation}\label{decom1}
h_{n+1}:=f-\sum_{k=1}^{n}h_{k}
\end{equation}
has no discontinuities on $\partial\Di$, that is, $h_{n+1}\in A_{0}$. Therefore $f\in A_{F}$.

Next, if $F$ consists of a single point, say $z_{0}$, then for a given $f\in AP(F)\cap H^{\infty}$ using the above argument we can find a function $h\in A_{\{z_{0},z_{1}\}}$ with a fixed $z_{1}\in\partial\Di$ such that $f-h$ is continuous on $F$. Let $g\in A_{0}$ be a function equal to $1$ at $F$ and $0$ at $z_{1}$.
Then $f-gh\in A_{0}$. This completes the first part of the proof.

Show now that $A_{F}\subset AP(F)\cap H^{\infty}$. Again, assume first that $F$ contains at least two points.
Let $e^{\lambda f}\in A_{F}$, $\lambda\in\Re$, where
${\rm Re}\ \!f$ is the characteristic function of an arc, say $[x,y]$, with $x,y\in F$. Let $\psi:\partial\Di\to\partial\Di$ be the restriction to $\partial\Di$ of a M\"{o}bius transformation sending $1$ to $x$ and $-1$ to $y$. Then by the definition 
$$
(f\circ\psi\circ\phi_{1}^{-1})(z)=-\frac{i}{\pi}{\rm Log}\ \!z +C,\ \ \ z\in\H_{+},
$$
for some constant $C$. Thus we have
$$
e^{(\lambda f\circ\psi\circ ({\rm Log}\circ\psi_{1})^{-1})(z)}=e^{\lambda C}e^{-i\lambda z/\pi},\ \ \
z\in\Sigma.
$$
This means that $e^{\lambda f\circ\psi}\in AP(\{-1,1\})\cap H^{\infty}$. Then by Corollary \ref{cor1} we get $e^{\lambda f}\in AP(F)\cap H^{\infty}$. Since $A_{F}$ is generated by $A_{0}$ and such functions $e^{\lambda f}$, we obtain the required implication.

If $F$ is a single point, then we must show that $ge^{\lambda f}\in AP(F)\cap H^{\infty}$, $\lambda\in\Re$, where ${\rm Re}\ \!f$ is the characteristic function of an arc with an endpoint at $F$ and $g\in A_{0}$ is such that $ge^{f}$ has discontinuity at $F$ only. The result follows easily from the previous part of the proof, because $e^{\lambda f}$ is almost periodic on $\partial\Di\setminus\{F,y\}$ for some $y$ and is continuous at $y$. 

This completes the proof of the lemma.\ \ \ \ \ $\Box$

As it was mentioned above the required statement of the theorem follows from Lemmas \ref{finite} and \ref{fin1}\ \ \ \ \  $\Box$\\
\\
{\bf 4.3. Proof of Example \ref{notequal}.} Using the bilinear transformation
$\phi_{1}:\Di\to\H_{+}$, see (\ref{fi}), that maps $1$ to $0\in\Re$ and $-1$ to $\infty$ we can transfer the problem to a similar one for functions on $\H_{+}$. Namely, let $\{x_{k}\}_{k\in\N}\subset\Re_{+}$ be the sequence converging to $0$, which is the image of the sequence $\{e^{it_{k}}\}_{k\in\N}\subset\partial\Di$ of the example under $\phi_{1}$. Let $H:\Re\to\{0,1\}$ be the Heaviside function (i.e., the characteristic function of $[0,\infty)$). Then the pullback by $\phi_{1}^{-1}$ of the function $u$ of the example to the boundary $\Re$ of $\H_{+}$ is the function
$$
\widetilde u(x):=\sum_{k=1}^{\infty}\alpha_{k}H(x-x_{k}),\ \ \ x\in\Re.
$$
We extend $\widetilde u$ to a harmonic function on $\H_{+}$ by the Poisson integral. Let $\widetilde v$ be the harmonic conjugate to the  extended function determined on $\Re$ by the formula
\begin{equation}\label{ex17}
\widetilde v(x)=\sum_{k=1}^{\infty}\alpha_{k}\frac{\ln|x-x_{k}|}{\pi},\ \ \ x\in\Re.
\end{equation}
We set $\widetilde h:=\widetilde u+i\widetilde v$. Then $\widetilde h$ is the pullback by $\phi_{1}^{-1}$ of a holomorphic function $h$ on $\Di$ such that ${\rm Re}\ \!h|_{\partial\Di}=u$.
Assume, to the contrary, that
$e^{h}\in A_{S}$. Since $e^{u}\in R_{S}$, this assumption implies that $e^{iv}\in AP(S)$, where $v=\phi_{1}^{*}(\widetilde v|_{\Re})$. Then according to the definition of the topology on ${\cal M}(AP(S))$, see section 3.3, the functions $\cos(\widetilde v(e^{t}))$ and $\sin(\widetilde v(e^{t}))$, $t\in\Re$, admit continuous extensions to $b\Re$
determined as follows.

If $\{s_{\alpha}\}\subset\Re$ is a net converging in $b\Re$ to a point $\eta\in b\Re$, then the values at $\eta$ of the extended functions are 
$$
\lim_{\alpha}\cos(\widetilde v(e^{s_{\alpha}}))\ \ \ {\rm and}\ \ \   \lim_{\alpha}\sin(\widetilde v(e^{s_{\alpha}})),
$$
respectively. In particular, this definition requires the existence of these limits.

Now, according to (\ref{ex17}) there are sequences of points
$\{x_{k}'\}_{k\in\N}$ and $\{x_{k}''\}_{k\in\N}$ in $\Re_{+}$ such that $x_{k}'$ and $x_{k}''$ are sufficiently close to $x_{k}$ and 
$$
\lim_{k\to\infty}\left|\frac{x_{k}'}{x_{k}''}\right|=1,\ \ \
\cos(\widetilde v(x_{k}'))=0,\ \ \ \cos(\widetilde v(x_{k}''))=1,\ \ \ k\in\N.
$$
We set $t_{k}':=\ln x_{k}'$ and $t_{k}'':=\ln x_{k}''$, $k\in\N$.
Assume without loss of generality that $\{t_{k}'\}$ forms a net converging in the topology of $b\Re$ to a point $\xi\in b\Re$ (for otherwise we replace $\{t_{k}\}$ by a proper subset satisfying this property). Since by
the definition $\lim_{k\to\infty} |t_{k}'-t_{k}''|=0$ and almost periodic functions on $\Re$ are uniformly continuous, the family $\{t_{k}''\}$ forms  a net with the same indeces as for the net formed by $\{t_{k}'\}$ whose limit in $b\Re$ is $\xi$, as well. Hence, we must have
$$
\lim_{k}\cos(\widetilde v(e^{t_{k}'}))=\lim_{k}\cos(\widetilde v(e^{t_{k}''})),
$$
a contradiction. Therefore $e^{h}\not\in A_{S}$.

Let now $f\in A_{0}$ be such that $f(1)=0$. Then the function $(fe^{h})|_{\partial\Di}$ is continuous at $1$. Thus $(fe^{h})|_{\partial\Di}$ can be uniformly approximated by constant functions on open arcs containing $1$. The same is true for each point $e^{it_{k}}$. This implies immediately that $fe^{h}\in A_{S}$.\ \ \ \ \ $\Box$  
\sect{\hspace*{-1em}. Proofs of Theorems \ref{te4} and \ref{te5}.}
{\bf 5.1. Proof of Theorem \ref{te4}.} Let us recall that $A_{S}$ is the uniform closure of the algebra generated by all possible $A_{F}$ with finite $F\subset S$. Therefore the maximal ideal space ${\cal M}(A_{S})$ of $A_{S}$ is the inverse limit of the maximal ideal spaces ${\cal M}(A_{F})$ of  $A_{F}$. In particular, if we will prove that $\Di$ is dense in each ${\cal M}(A_{F})$, then by the definition of the inverse limit this will imply that $\Di$ is dense in ${\cal M}(A_{S})$, as required. Thus it suffices to prove the theorem for $A_{F}$ with $F=\{z_{1},\dots, z_{n}\}\subset \partial\Di$.  
\begin{Th}\label{corfin}
$\Di$ is dense in ${\cal M}(A_{F})$.
\end{Th}
{\bf Proof.}
Let ${\cal I}_{k}\subset A_{F}$ be the closed ideal consisting of functions that are continuous and equal to 0 at $z_{k}$. By $A_{k}$ we denote the quotient Banach algebra $A_{F}/{\cal I}_{k}$ equipped with the quotient norm. Let us recall that ${\cal M}_{z_{k}}$ is the maximal ideal space of the algebra ${\cal A}_{z_{k}}$ which is the uniform closure of the algebra of continuous functions on $\Di$ almost periodic in circular neighbourhoods of $z_{k}$. Also, according to Lemma \ref{le3} (b) there is a natural continuous projection $p_{z_{k}}:{\cal M}_{z_{k}}\to\overline\Di$ and $p_{z_{k}}^{-1}(z_{k})$ is homeomorphic to
${\cal M}(AP_{\cal O}(\Sigma))$. Moreover, by Lemma \ref{le3} (c) each $f\in A_{F}$ is extended to a continuous function on ${\cal M}_{z_{k}}$  holomorphic on $p_{z_{k}}^{-1}(z_{k})$. Hence, there is a continuous map $H_{k}:{\cal M}_{z_{k}}\to {\cal M}(A_{F})$ whose image coincides with the closure of
$\Di$. Moreover, according to the decomposition obtained in the proof of Lemma \ref{fin1}, see (\ref{decom1}), $H_{k}$ maps $p_{z_{k}}^{-1}(z_{k})$ homeomorphically onto its image.
\begin{Lm}\label{quot}
Let $\phi_{k}:A_{F}\to AP_{\cal O}(\Sigma)$ be the composition of the extension homomorphism of functions from $A_{F}$ to ${\cal M}_{z_{k}}$ and the restriction homomorphism of functions on ${\cal M}_{z_{k}}$ to
$p_{z_{k}}^{-1}(z_{k})$. Then $Ker\ \!\phi_{k}={\cal I}_{k}$ and $A_{k}$ is isomorphic to $AP_{\cal O}(\Sigma)$.
\end{Lm}
{\bf Proof.} Clearly, ${\cal I}_{k}\subset Ker\ \!\phi_{k}$. Let us check the converse implication. Let $f\in Ker\ \!\phi_{k}$.
As it follows from the proof of Lemma \ref{fin1}, see (\ref{decom1}), there are linear continuous operators $T_{k}:AP_{\cal O}(\Sigma)\to A_{\{z_{k}\}}\subset A_{F}$  such that $\phi_{k}\circ T_{k}=id$. Moreover, 
$T_{0}:=I-\sum_{k=1}^{n}T_{k}\circ\phi_{k}$, where $I$ is the identity map, maps $A_{F}$ onto $A_{0}$. In particular, we have 
$T_{k}(\phi_{k}(f))=0$. Thus $f=-T_{0}(f)+\sum_{s\neq k}T_{s}(\phi_{s}(f))$.
Since $\phi_{k}(f)=0$, this implies that $f$ is continuous and equal to $0$ at $z_{k}$. Now from the formula $\phi_{k}\circ T_{k}=id$ we obtain that
$A_{k}$ is isomorphic to $AP_{\cal O}(\Sigma)$.\ \ \ \ \ \ $\Box$

Let $i:A_{0}\hookrightarrow A_{F}$ be the natural inclusion. Its dual determines a continuous surjective map $a_{F}:{\cal M}(A_{F})\to\overline\Di$.
Next, taking the dual map to $\phi_{k}$ we obtain that each ${\cal M}(AP_{\cal O}(\Sigma))$ is embedded into ${\cal M}(A_{F})$, its image coincides with $H_{k}(p_{z_{k}}^{-1}(z_{k}))$ and $a_{F}$ maps $H_{k}(p_{z_{k}}^{-1}(z_{k}))$ to $z_{k}$.

Let $\xi\in {\cal M}(A_{F})$ and ${\mathfrak{m}}:={\rm Ker}~\xi\subset A_{F}$ be the corresponding maximal ideal. Assume first that 
there is $k$ such that  ${\cal I}_{k}\subset {\mathfrak m}$. Then ${\mathfrak m}_{k}=\phi_{k}({\mathfrak m})$ is a maximal ideal of $A_{k}$. By $\xi_{k}\in {\cal M}(AP_{\cal O}(\Sigma))$ we denote the character corresponding to ${\mathfrak m}_{k}$. Then $\xi=\phi^{*}(\xi_{k})\in H_{k}(p_{z_{k}}^{-1}(z_{k}))$. Now, by the definition of $H_{k}$ the point $\xi$ belongs to the closure of $\mathbb D$ in ${\cal M}(A_{F})$. We continue with the following lemma.
\begin{Lm}\label{inters}
Assume that a maximal ideal ${\mathfrak m}$ of $A_{F}$ does not contain any of ${\cal I}_{k}$. Then ${\mathfrak m}$ does not contain $\cap_{1\leq k\leq n}{\cal I}_{k}$, as well.
\end{Lm}
{\bf Proof.} Suppose, to the contrary, that $\cap_{1\leq k\leq n}{\cal I}_{k}\subset {\mathfrak m}$. Let $x_{k}\in {\cal I}_{k}$, $1\leq k\leq n$, be such that $x_{k}\not\in {\mathfrak m}$. Since ${\cal I}_{k}$ are ideals, $x_{1}\cdots x_{n}\in\cap_{1\leq k\leq n}{\cal I}_{k}$.
Thus $x_{1}\cdots x_{n}\in {\mathfrak m}$. Since ${\mathfrak m}$ is a prime ideal, there is some $k$ so that $x_{k}\in {\mathfrak m}$, a contradiction.\ \ \ \ \ $\Box$

From this lemma it follows that in order to prove the theorem it remains to consider the case ${\mathfrak m}\not\in\cap_{1\leq k\leq n}{\cal I}_{k}$.
Observe that $\cap_{1\leq k\leq n}{\cal I}_{k}$ consists of all functions from $A_{0}$ that vanish on $F$.  Thus there is $f\in \cap_{1\leq k\leq n}{\cal I}_{k}$ such that
$f(\xi)\neq 0$. This implies that $a_{F}(\xi)\not\in  F$. For every
$g\in A_{F}$, let us consider the function $gf$. By the definition $gf\in A_{0}$. Thus we have 
$$
g(a_{F}(\xi))f(a_{F}(\xi))=(gf)(a_{F}(\xi))=(gf)(\xi)=g(\xi)f(\xi)=g(\xi)f(a_{F}(\xi)).
$$
Equivalently,
$$
g(\xi)=g(a_{F}(\xi))\ \ \ {\rm for\ all}\ \ \ g\in A_{F}.
$$
This implies that $a_{F}^{-1}(a_{F}(\xi))=\{\xi\}$. Therefore  $a_{F}:{\cal M}(A_{F})\setminus a_{F}^{-1}(F)\to\overline \mathbb D\setminus F$ is a homeomorphism. In particular, $\xi$ belongs to the closure of $\mathbb D$.

The proof of Theorem \ref{te4} is complete.\ \ \ \ \ $\Box$\\
\\
{\bf 5.2. Proof of Theorem \ref{te5}.} 
Statements (1) and (3) of the theorem follow easily from similar statements for $a_{F}$ with a finite subset $F\subset S$ proved in section 5.1 and the properties of the inverse limit. Let us prove (3). 
We first prove the statement for a finite subset $F\subset S$.

Since $A_{0}\subset A_{F}$ and each function from $A_{0}$ attains the maximum of modulus on $\partial\mathbb D$, $K_{F}\subset a_{F}^{-1}(\partial\mathbb D)$. As it follows from Theorem \ref{te3} $A_{F}\hookrightarrow AP(F)$. Also, the extensions of functions from $A_{F}$  to ${\cal M}(AP(F))$ separate points there. Therefore ${\cal M}(AP(F))$ is embedded into ${\cal M}(A_{F})$.  Identifying ${\cal M}(AP(F))$ with its image under this embedding we have
${\cal M}(AP(F))\subset a_{F}^{-1}(\partial\mathbb D)$. Since $a_{F}^{-1}(\partial\mathbb D)\setminus a_{F}^{-1}(F)\to\partial\mathbb D\setminus F$ is a homeomorphism and each $z\in\partial\mathbb D$ is a peak point for $A_{0}$, the set $K_{F}$ contains the closure of $\mathbb D\setminus F$ which is, by definition, coincides with ${\cal M}(AP(F))$. Assume that there is 
$\xi\in K_{F}\setminus {\cal M}(AP(F))$. Then  $a_{F}(\xi):=z^{*}\in F$. Further, identifying $a_{S}^{-1}(z^{*})$ with ${\cal M}(AP_{\cal O}(\Sigma))$
we get from our assumption that $\xi\in i_{\eta}(\Sigma_{0})$ for some $\eta\in b\mathbb Z$, see section 2; here $\Sigma_{0}$ is the interior of $\Sigma$. Then, because $i_{\eta}(\Sigma_{0})$ is dense in ${\cal M}(AP_{\cal O}(\Sigma))$, by the maximum modulus principle each function $f\in A_{F}$
for which $\max_{\mathbb D}|f|=|f(\xi)|$ is constant on $a_{S}^{-1}(z^{*})$. Thus
such $f$ admits the maximum of modulus on ${\cal M}(AP(F))$ also. This contradicts to the minimality of $K_{F}$. Therefore $K_{F}= {\cal M}(AP(F))$. 

Further, ${\cal M}(AP(S))$ is the inverse limit of compact sets ${\cal M}(AP(F))$ for all finite $F\subset S$. As before we naturally identify ${\cal M}(AP(S))$
with a subset of ${\cal M}(A_{S})$. Then since $A_{S}$ is the uniform closure of the algebra generated by all possible $A_{F}$ with finite $F\subset S$, by the definition of the inverse limit $K_{S}\subseteq {\cal M}(AP(S))$.  But in fact $K_{S}$ coincides with ${\cal M}(AP(S))$ because otherwise its projection to some of ${\cal M}(A_{F})$ is a boundary of
$A_{F}$ and a proper subset of ${\cal M}(AP(F))$, a contradiction.
\ \ \ \ \ $\Box$


\end{document}